\documentclass[a4paper,11pt]{article}
\usepackage{amsmath,amsthm,amssymb,amsfonts,fullpage}

\newtheorem{thm}{Theorem}[section]

\newtheorem{assu}{Assumption}
\newtheorem{cor}[thm]{Corollary}

\newtheorem{lem}[thm]{Lemma}
\newtheorem{propn}[thm]{Proposition}
\newtheorem{rem}[thm]{Remark}

\newcommand{\Ecal}   {\mathcal{E}}
\newcommand{\Fcal}   {\mathcal{F}}

\newcommand{\bR}{{\mathbb R}}

\def\eq#1\en{\begin{equation}#1\end{equation}}

\def\Reff{R_{\rm eff}}

\def\eps{\varepsilon}
\def\lam{\lambda}

\newcommand{\al}{\alpha}

\renewcommand{\to}      {\rightarrow}

\def\qed{{\hfill $\square$ \bigskip}}
\makeatletter
\@addtoreset{equation}{section}

\makeatother

\begin{document}

\title{Convergence of mixing times for sequences of \\ random walks on finite graphs}
\author{D.A. Croydon\footnote{Dept of Statistics, University of Warwick, Coventry, CV4 7AL, United Kingdom; d.a.croydon@warwick.ac.uk.}, ~B.M. Hambly\footnote{Mathematical Institute, 24-29 St Giles', Oxford, OX1
3LB, United Kingdom; hambly@maths.ox.ac.uk.}~ and ~T. Kumagai\footnote{RIMS, Kyoto University, Kyoto 606-8502, Japan; kumagai@kurims.kyoto-u.ac.jp.}}
\maketitle

\begin{abstract}
We establish conditions on sequences of graphs which ensure that the mixing times of the random walks on the
graphs in the sequence converge. The main assumption is that the graphs, associated measures and heat kernels
converge in a suitable Gromov-Hausdorff sense. With this result we are able to establish the convergence of the
mixing times on the largest component of the Erd\H{o}s-R\'enyi random graph in the critical window, sharpening
previous results for this random graph model. Our results also enable us to establish convergence in a number
of other examples, such as finitely ramified fractal graphs, Galton-Watson trees and the range of a
high-dimensional random walk.
\end{abstract}

\section{Introduction}

The geometric and analytic properties of random graphs have been the subject of much recent research.
One strand of this development has been to examine sequences of random subgraphs of vertex transitive
graphs that are, in some sense, at or near criticality. A key example is the percolation model and, for bond
percolation above the upper critical dimension, we expect to see mean-field behavior in the sequence of
finite graphs in the critical window. That is, the natural scaling exponents for the volume and diameter of
the graph and for the mixing time are of the same order as those for the Erd\H{o}s-R\'enyi random graph
in the critical window, as given in \cite{NacPer08}.

This mean-field behavior is seen in other natural models of sequences of critical random graphs. For
example \cite{BCHSS} obtained general conditions for the geometric properties of percolation clusters
on sequences of finite graphs and discussed examples such as the high dimensional torus and the
$n$-cube, while the random walk on critical percolation clusters on the high-dimensional torus is treated in \cite{HeyvdH09}.
Motivated by these results we will focus on the asymptotic behavior of mixing times for random walks
on sequences of finite graphs. We consider general sequences of graphs but under some strong conditions
which will enable us to establish the convergence of the mixing time.

In order to demonstrate our main result we consider the Erd\H{o}s-R\'enyi random graph. Let  $G(N,p)$
be the random subgraph of the complete graph on $N$ labeled vertices $\{1,\dots, N\}$ in which each edge is
present with probability $p$ independently of the other edges.  It is a classical result that if
we set $p=c/N$, then as $N\to\infty$, if $c>1$ there
is a giant component containing a positive fraction of the vertices, while for $c<1$ the largest component is
of size $\log{N}$. However, if $p=N^{-1}+\lambda N^{-4/3}$ for
some $\lambda\in \mathbb{R}$, we are in the so-called critical window, and it is known that the largest
connected component $\mathcal{C}^N$, is of order $N^{2/3}$. The recent work of \cite{ABG} has shown that
the scaling limit of the graph, $\mathcal{M}$, exists and can be constructed from the continuum random tree.

For the Erd\H{o}s-R\'enyi random graph above criticality, \cite{FR} and \cite{BKW} established mixing time
bounds for the simple random walk on the giant component. The simple random walk on this graph is
the discrete time Markov chain with transition probabilities determined by $p(x,y) = 1/$deg$(x)$ for all $y$
such that $(x,y)$ is an edge in $\mathcal{C}^N$. For the random graph in the critical window, the
following result on the mixing time $t_{\rm mix}^1(\mathcal{C}^N)$  (a precise
definition will be given later in (\ref{mixdefdis}), see also Remark~1.3) of the lazy random walk (a version of the
simple random walk which remains at its current vertex with probability 1/2, otherwise
it moves as the simple random walk) was obtained by Nachmias and Peres (\cite[Theorem 1.1]{NacPer08}).
{\thm
Let ${\cal C}^N$ be the largest connected
component of $G(N,(1+\lambda N^{-1/3})/N)$ for some $\lambda\in {\mathbb R}$.
Then, for any $\epsilon>0$, there exists $A=A(\epsilon,\lambda)<\infty$ such that
for all large $N$,
\[ P(t_{\rm mix}^1({\cal C}^N) \notin [A^{-1}N,A N])<\epsilon.\]}

It is natural to ask for more refined results on the behavior of the family of mixing times. The purpose of this paper is
to give a general criteria for the convergence of mixing times for a sequence of simple random walks on finite graphs in the
setting where the graphs can be embedded nicely in a compact metric space. Due to the recent work of \cite{ABG}
and \cite{Croydoncrg} we can apply our main result to the case of the Erd\H{o}s-R\'enyi random graph, to obtain the
following result.

\begin{thm}
Fix $p\in [1,\infty]$. If $t_{\rm mix}^{p}(\rho^N)$ is the $L^p$-mixing time of the simple random walk on
$\mathcal{C}^N$ started from its root $\rho^N$, then
\[N^{-1}t_{\rm mix}^{p}(\rho^N) {\rightarrow} t_{\rm mix}^{p}(\rho),\]
in distribution, where the random variable $t_{\rm mix}^{p}(\rho)\in(0,\infty)$ is the  $L^p$-mixing
time of the Brownian motion on $\mathcal{M}$ started from $\rho$.
\end{thm}

We will later illustrate our main result with a number of other examples of random walks on sequences of finite graphs.
In order to state it, though, we start by describing the general framework in which we work.
Firstly, let $(F,d_F)$ be a compact metric space
and let $\pi$ be a non-atomic Borel probability measure on $F$ with full support. We will assume that
balls $B_F(x,r):=\{y\in F:d_F(x,y)<r\}$ are $\pi$-continuity sets (i.e. $\pi(\partial B_F(x,r))=0$ for every
$x\in F$, $r>0$). Secondly, take $X^F=(X^F_t)_{t\geq 0}$ to be a $\pi$-symmetric Hunt process
on $F$ (for definition and properties see \cite{FOT}), which will typically be the Brownian motion on the
limit of the sequence of graphs. We suppose the following:
\begin{eqnarray}
&&\mbox{\textbullet\hspace{5pt} $X^F$ is conservative, i.e. its semigroup $(P_t)_{t\geq 0}$ satisfies $P_t1=1$, $\pi$-a.e., $\forall t>0$,}\label{asum-a}\\
&&\mbox{\textbullet\hspace{5pt} there exists a jointly continuous transition density $(q_t(x,y))_{x,y\in F,t>0}$
 of $X^F$,}\hspace{30pt}\label{asum-b}\\
&&\mbox{\textbullet\hspace{5pt} for every $x,y\in F$ and $t>0$, $q_t(x,y)>0$,}\label{asum-c}\\
&&\mbox{\textbullet\hspace{5pt} for every $x\in F$ and $t>0$, $q_t(x,\cdot)$ is not identically equal to 1,}\label{asum-d}
\end{eqnarray}
where conditions \eqref{asum-c} and \eqref{asum-d} are assumed to exclude various trivial cases, and by
transition density we mean the kernel $q_t(x,y)$ such that
\[{\mathbf E}_x[f(X^F_t)]=\int_Fq_t(x,y)f(y)\pi(dy),\]
for all bounded continuous function $f$ on $F$. Furthermore, we will say that the transition density
$(q_t(x,y))_{x,y\in F,t>0}$ converges to stationarity in an $L^p$ sense for some $p\in[1,\infty]$ if it holds that
\begin{equation}\label{erg}
\lim_{t\rightarrow\infty}D_p(x,t)=0,
\end{equation}
for every $x\in F$, where $D_p(x,t):=\|q_t(x,\cdot)-1\|_{L^p(\pi)}$. If this previous condition is satisfied, then it is
possible to check that the $L^p$-mixing time of $F$,
\begin{equation}\label{mixingdef}
t_{\rm mix}^p(F):=\inf\left\{t>0:\sup_{x\in F}D_p(x,t)\leq 1/4\right\},
\end{equation}
is a finite quantity (see Section \ref{mixingtime}). Finally, note that $t_{\rm mix}^p(F)\le t_{\rm mix}^{p'}(F)$ for
$p\le p'$, which can easily be shown using the H\"older inequality.

We continue by introducing some general notation for graphs and their associated random walks. First,
fix $G=(V(G),E(G))$ to be a finite connected graph with at least two vertices, where $V(G)$ denotes the
vertex set and $E(G)$ the edge set of $G$, and suppose $d_G$ is a metric on $V(G)$. In some examples,
$d_G$ will be a rescaled version of the usual shortest path graph distance, by which we mean that
$d_G(x,y)$ is some multiple of the number of edges in the shortest path from $x$ to $y$ in $G$, but
this is not always the most convenient choice. Define a symmetric weight function $\mu^G:V(G)^2
\rightarrow \mathbb{R}_+$ that satisfies $\mu^G_{xy}>0$ if and only if $\{x,y\}\in E(G)$. The discrete
time random walk on the weighted graph $G$ is then the Markov chain $((X^G_m)_{m\geq 0}, \mathbf{P}^G_x,x
\in V(G))$ with transition probabilities $(P_G(x,y))_{x,y\in V(G)}$ defined by $P_G(x,y):={\mu^G_{xy}}/{\mu^G_{x}}$,
where $\mu^G_x:=\sum_{y\in V(G)}\mu^G_{xy}$. If we define a measure $\pi^G$ on $V(G)$ by setting, for
$A\subseteq V(G)$, $\pi^G(A):=\sum_{x\in A}\mu^G_x/\sum_{x\in V(G)}\mu^G_x$,
then $\pi^G$ is the invariant probability measure for $X^G$. The transition density of $X^G$, with respect to
$\pi^G$, is given by $(p^G_m(x,y))_{x,y\in V(G),m\geq 0}$, where
\[p_m^G(x,y):=\frac{\mathbf{P}^{G}_x(X_m=y)}{\pi^G(\{y\})}.\]
Due to parity concerns for bipartite graphs, we will consider a smoothed version of this function
$(q^G_m(x,y))_{x,y\in V(G),m\geq 0}$ obtained by setting
\begin{equation}\label{smooth}
q^G_m(x,y):=\frac{p^G_m(x,y)+p^G_{m+1}(x,y)}{2},
\end{equation}
and define the $L^p$-mixing time of $G$ by
\begin{equation}\label{mixdefdis}
t_{\rm mix}^{p}(G):=\inf\left\{m> 0:\sup_{x\in V(G)}D^{G}_p(x,m)\leq 1/4\right\},\end{equation}
where $D^{G}_p(x,m):=\|q^G_m(x,\cdot)-1\|_{L^p(\pi^G)}$. Finally, in the case that we are considering a sequence
of graphs $(G^N)_{N\geq 1}$, we will usually abbreviate $\pi^{G^N}$ to $\pi^N$ and $q^{G^N}$ to $q^N$, etc.

{\rem {\rm
In \cite{NacPer08}, the mixing time of $\mathcal{C}^N$ is defined in terms of the total variation distance, that is
\begin{equation}\label{N+Pnomix}
T_{\rm mix}({\cal C}^N)=\min\{t:\|P_t(x,\cdot)-\pi(\cdot)\|_{\rm TV}\le 1/8,~~\forall x\in V({\cal C}^N)\},\end{equation}
where $P_t(x, B)=\sum_{y\in B}p_t^N(x,y)\pi(y)$ for $B\subset V({\cal C}^N)$, $p_t^N(x,y)$ is the transition density
for the random walk and $\|\mu-\nu\|_{\rm TV}=\max_{B\subset V({\cal C}^N)}|\mu(B)-\nu(B)|$
for probability measures $\mu,\nu$ on $V({\cal C}^N)$. (To be precise, $1/8$ in \eqref{N+Pnomix} is $1/4$ in
\cite{NacPer08}, but this only affects the constants in the results.)  However, noting that
\[\|\mu-\nu\|_{\rm TV}=\frac 12\sum_{x\in V({\cal C}^N)}|\mu(\{x\})-\nu(\{x\})|,\]
(see, for example \cite[Proposition 4.2]{LPW}), one sees that $T_{\rm mix}({\cal C}^N)=t_{\rm mix}^{1}({\cal C}^N)$. Also note that \cite{NacPer08} considers the lazy walk on the graph to avoid parity issues, but the same techniques will apply to the mixing time defined in terms of the smoothed heat kernel introduced at (\ref{smooth}).}}

\bigskip

We are now ready to state the assumption under which we are able to prove the convergence of mixing times
for the random walks on a sequence of graphs. This captures the idea that, when suitably rescaled, the discrete state
spaces, invariant measures and transition densities of a sequence of graphs converge to $(F,d_F)$, $\pi$ and
$(q_t(x,y))_{x,y\in F,t>0}$, respectively. Its formulation involves a spectral Gromov-Hausdorff topology, the
definition of which is postponed until Section \ref{gh}, and a useful sufficient condition for it will be given in
Proposition~\ref{local} below. Note that we extend the definition of the discrete transition densities on graphs to all
positive times by linear interpolation of $(q^G_m(x,y))_{m\geq 0}$ for each pair of vertices $x,y\in V(G)$.
Note also that the extended transition densities are different from those of continuous time Markov chains.

{\assu\label{assu1} $(G^N)_{N\geq 1}$ is a sequence of finite connected graphs with at least two vertices
for which there exists a sequence  $(\gamma(N))_{N\geq 1}$ such that, for any compact interval $I\subset (0,\infty)$,
\[{\left(\left(V(G^N),d_{G^N}\right), \pi^N, \left(q^N_{\gamma(N)t}(x,y)\right)_{x,y\in V(G^N),t\in I}\right)}
\rightarrow\left(\left(F,d_F\right),\pi,(q_t(x,y))_{x,y\in F,t\in I}\right)\]
in a spectral Gromov-Hausdorff sense.}
\bigskip

In the case where we have random graphs, we will typically assume that we have the above convergence holding in distribution. Our main conclusion is then the following.

{\thm \label{main} Suppose that Assumption \ref{assu1} is satisfied. If $p\in[1,\infty]$ is such that the transition
density $(q_t(x,y))_{x,y\in F,t>0}$ converges to stationarity in an $L^p$ sense, then $t_{\rm mix}^{p}(F)\in (0,\infty)$
and
\begin{equation}\label{mainconv}
\gamma(N)^{-1}t_{\rm mix}^{p}(G^N)\rightarrow t_{\rm mix}^{p}(F).
\end{equation}}

In Section \ref{fixed}, we will explain how to derive a variation of Theorem \ref{main}
that concerns the convergence of mixing times of processes started at a distinguished point in the state space.

We emphasize that a key part of our paper is to verify Assumption \ref{assu1} and apply Theorem \ref{main}
in various interesting examples (including the Erd\H{o}s-R\'enyi random graphs in the critical window as mentioned
above). Therefore, we devote considerable space to applying our results to such examples.

The organization of the paper is as follows. In Section 2, we give a precise definition of the spectral Gromov-Hausdorff
convergence and give some of its basic properties. In Section 3, we prove Theorem~\ref{main} and derive a variation
of the theorem for distinguished starting points. Some sufficient conditions for \eqref{asum-a}-\eqref{erg} are
given in Section 4. A selection of examples where the assumptions of Theorem~\ref{main} can be verified, and
hence we have convergence of the mixing time sequence, are given in Section 5. In Section 6 we introduce some
geometric conditions on graphs for upper and lower bounds on the mixing times for the
corresponding symmetric Markov chains. We use these ideas to derive tail estimates of mixing times on random graphs
in the case of the continuum random tree and the Erd\H{o}s-R\'enyi random graph. The proofs of these results
can be found in the Appendix.

\section{Spectral Gromov-Hausdorff convergence}\label{gh}

The aim of this section is to define a spectral Gromov-Hausdorff distance on triples consisting of
a metric space, a measure and a heat kernel-type function that will allow us to make
Assumption \ref{assu1} precise. We will also derive an equivalent characterization of this assumption
that will be applied in the subsequent section when proving our mixing time convergence result,
and present a sufficient condition for Assumption \ref{assu1} that will be useful when it comes to
checking it in examples. Note that we do not need to assume  \eqref{asum-c}, \eqref{asum-d} in
this section, and only use \eqref{asum-a} to deduce Proposition \ref{local} from a result of
\cite{CHLLT}.

First, for a compact interval $I\subset(0,\infty)$, let $\tilde{\mathcal{M}}_I$ be the collection of
triples of the form $(F,\pi,q)$, where $F=(F,d_F)$ is a non-empty compact metric space, $\pi$ is
a Borel probability measure on $F$ and $q=(q_t(x,y))_{x,y\in F, t\in I}$ is a jointly continuous
real-valued function of $(t,x,y)$. We say two elements, $(F,\pi,q)$ and $(F',\pi',q')$, of
$\tilde{\mathcal{M}}_I$ are equivalent if there exists an isometry $f:F\rightarrow {F'}$ such that
$\pi\circ f^{-1}=\pi'$ and $q_t'\circ f=q_t$ for every $t\in I$, by which we mean $q_t'(f(x),f(y))=q_t(x,y)$
for every $x,y\in F$, $t\in I$. Define $\mathcal{M}_I$ to be the set of equivalence classes of
$\tilde{\mathcal{M}}_I$ under this relation. We will often abuse notation and identify an equivalence
class in ${\mathcal{M}}_I$ with a particular element of it. Now, set
\begin{eqnarray*}
\lefteqn{\Delta_I\left((F,\pi,q),(F',\pi',q')\right)}\\
&:=&\inf_{Z,\phi,\phi',\mathcal{C}}\left\{d_H^Z(\phi(F),\phi'(F'))+d_P^Z(\pi\circ\phi^{-1},\pi'\circ\phi'^{-1})
\vphantom{+\sup_{(x,x'),(y,y')\in\mathcal{C}}\left(d_Z(\phi(x),\phi'(x'))+d_Z(\phi(y),\phi'(y'))+\left|q_t(x,y)-q't(x',y')\right|\right)}\right.\\
&&\left.+\sup_{(x,x'),(y,y')\in\mathcal{C}}\left(d_Z(\phi(x),\phi'(x'))+d_Z(\phi(y),\phi'(y'))+\sup_{t\in I}\left|q_t(x,y)-q'_t(x',y')\right|\right)\right\},
\end{eqnarray*}
where the infimum is taken over all metric spaces $Z=(Z,d_Z)$, isometric embeddings $\phi:F\rightarrow Z$, $\phi':F'\rightarrow Z$, and correspondences $\mathcal{C}$ between $F$ and $F'$, $d_H^Z$ is the Hausdorff distance between compact subsets of $Z$, and $d_P^Z$ is the Prohorov distance between Borel probability measures on $Z$. Note that, by a correspondence $\mathcal{C}$ between $F$ and ${F'}$, we mean a subset of $F\times {F'}$ such that for every $x\in F$ there exists at least one $x'\in {F'}$ such that $(x,x')\in\mathcal{C}$ and conversely for every $x'\in {F'}$ there exists at least one $x\in F$ such that $(x,x')\in \mathcal{C}$.

In the following lemma, we check that the above definition gives us a metric and that the corresponding space is separable. (The latter fact will be useful when it comes to making convergence in distribution statements regarding the mixing times of sequences of random graphs, as is done in Sections \ref{treesec} and \ref{ersec}, for example). Before this, however, let us make a few remarks about the inspiration for the distance in question. In the infimum characterizing $\Delta_I$, the first term is simply that used in the standard Gromov-Hausdorff distance  (see \cite[Definition 7.3.10]{BBI}, for example). In fact, as far as the topology is considered, this term could have been omitted since it is absorbed by the other terms in the expression, but we find that it is technically convenient and somewhat instructive to maintain it. The second term is that considered by the authors of \cite{GPW} in defining their `Gromov-Prohorov' distance between metric measure spaces. The final term is closely related to one used in
\cite[Section 6]{LegallDuquesne}
when defining a distance between spatial trees -- real trees equipped with a continuous function. Indeed, the notion of a correspondence is quite standard in the Gromov-Hausdorff setting as a way to relate two compact metric spaces. One can, for example, alternatively define the Gromov-Hausdorff distance between compact metric spaces as half the infimum of the distortion of the correspondences between them (see \cite[Theorem 7.3.25]{BBI}).

{\lem \label{ghlem} For any compact interval $I\subset (0,\infty)$, $(\mathcal{M}_I,\Delta_I)$ is a separable metric space.}
\begin{proof} Fix a compact interval $I\subset (0,\infty)$. That $\Delta_I$ is a non-negative function and is symmetric is obvious. To prove that it is also the case that $\Delta_I\left((F,\pi,q),(F',\pi',q')\right)<\infty$ for any choice of $(F,\pi,q),(F',\pi',q')\in\mathcal{M}_I$, simply consider $Z$ to be the disjoint union of $F$ and $F'$, setting $d_Z(x,x'):={\rm diam}(F,d_F)+{\rm diam}(F',d_F')$ for any $x\in F,x'\in F'$, and suppose that $\mathcal{C}=F\times F'$.

We next show that $\Delta_I$ is positive definite. Suppose $(F,\pi,q),(F',\pi',q')\in\mathcal{M}_I$ are such that $\Delta_I\left((F,\pi,q),(F',\pi',q')\right)=0$. For every $\varepsilon>0$, we can thus choose $Z,\phi,\phi',\mathcal{C}$ such that the sum of quantities in the defining infimum of $\Delta_I$ is bounded above by $\varepsilon$. Moreover, there exists a $\delta\in(0,\varepsilon]$ such that
\begin{equation}\label{qcont}
\sup_{\substack{x_1,x_2,y_1,y_2\in F:\\d_F(x_1,x_2),d_F(y_1,y_2)\leq \delta}}\sup_{t\in I}\left|q_t(x_1,y_1)-q_t(x_2,y_2)\right|\leq \varepsilon.
\end{equation}
Now, let $(x_i)_{i=1}^{\infty}$ be a dense sequence of disjoint elements of $F$ (in the case $F$ is finite, we suppose that the sequence terminates after having listed all of the elements of $F$). By the compactness of $F$, there exists an integer $N_\varepsilon$ such that $(B_F(x_i,\delta))_{i=1}^{N_{\varepsilon}}$ is a cover for $F$. Define $A_1:=B_F(x_1,\delta)$, and $A_i:=B_F(x_i,\delta)\backslash\cup_{j=1}^{i-1}B_F(x_i,\delta)$ for $i=2,\dots,N_\varepsilon$, so that $(A_i)_{i=1}^{N_\varepsilon}$ is a disjoint cover of $F$, and then consider a function $f_\varepsilon:F\rightarrow F'$ obtained by setting
\[f_\varepsilon(x):=x_i'\]
on $A_i$, where $x_i'$ is chosen such that $(x_i,x_i')\in\mathcal{C}$ for each $i=1,\dots,N_\varepsilon$. Clearly, by definition, $f_\varepsilon$ is a measurable function. It is further the case that it satisfies, for any $x\in F$,
\[d_Z(\phi(x),\phi'(f_\varepsilon(x)))\leq d_Z(\phi(x),\phi(x_i))+d_Z(\phi(x_i),\phi'(x_i'))\leq 2\varepsilon,\]
where, in the above, we assume that $i\in \{1,\dots,N_\varepsilon\}$ is such that $x\in A_i$. From this, it readily follows that:
\begin{equation}\label{fact1}
\sup_{x,y\in F}\left|d_F(x,y)-d_{F'}(f_\varepsilon(x),f_\varepsilon(y))\right|\leq 4\varepsilon
\end{equation}
and
\begin{equation}\label{fact2}
d_P^{F'} (\pi\circ f^{-1}_{\varepsilon},\pi')\leq 3\varepsilon,
\end{equation}
where $d_P^{F'}$ is the Prohorov distance on $F'$. By applying (\ref{qcont}), we also have that
\begin{equation}\label{fact3}
\sup_{x,y\in F,t\in I}\left|q_t(x,y)-q'_t(f_\varepsilon(x),f_\varepsilon(y))\right|\leq 2\varepsilon.
\end{equation}
To continue, we use a diagonalization argument to deduce the existence of a sequence $(\varepsilon_{n})_{n\geq 1}$ such that $f_{\varepsilon_n}(x_i)$ converges to some limit $f(x_i)\in F'$ for every $i\geq 1$. From (\ref{fact1}), we obtain that $d_{F'}(f(x_i),f(x_j))=d_F(x_i,x_j)$ for every $i,j\geq 1$, and so we can extend the map $f$ continuously to the whole of $F$ (\cite[Proposition 1.5.9]{BBI}). This construction immediately implies that $f$ is distance preserving. Moreover, reversing the roles of $F$ and $F'$, we are able to find a distance preserving map from $F'$ to $F$. Hence $f$ must be an isometry. To check that  $(F,\pi,q)$ and $(F',\pi',q')$ are equivalent, it therefore remains to check that $\pi\circ f^{-1}=\pi'$ and  $q_t'\circ f=q_t$ for every $t\in I$. Fix $\varepsilon>0$ and recall that the definition of  $(x_i)_{i=1}^{N_\varepsilon}$ means that it is an $\varepsilon$-net for $F$. Let $\varepsilon'\in(0, \varepsilon]$ be such that $d_{F'}(f_{\varepsilon'}(x_i),f(x_i))\leq \varepsilon$ for every $i=1,\dots,N_{\varepsilon}$. Then,
\begin{equation}\label{fact4}
d_{F'}(f_{\varepsilon'}(x),f(x))\leq d_{F'}(f_{\varepsilon'}(x),f_{\varepsilon'}(x_i))+
d_{F'}(f_{\varepsilon'}(x_i), f(x_i))+d_{F'}( f(x_i),f(x))
\leq 7\varepsilon,
\end{equation}
where we are again assuming that $i\in \{1,\dots,N_\varepsilon\}$ is such that $x\in A_i$, and have applied (\ref{fact1}) and the distance-preserving property of $f$. In particular, this implies that
\[d_P^{F'} (\pi\circ f^{-1},\pi')\leq d_P^{F'} (\pi\circ f^{-1},\pi\circ f_{\varepsilon'}^{-1})+ d_P^{F'} (\pi\circ f^{-1}_{\varepsilon'},\pi')\leq 10\varepsilon,\]
where we use (\ref{fact2}) to deduce the second inequality. Since $\varepsilon>0$ was arbitrary, this yields that $\pi\circ f^{-1}=\pi'$. Finally, (\ref{fact3}) and (\ref{fact4}) imply that
\[\sup_{x,y\in F,t\in I}\left|q_t(x,y)-q'_t(f(x),f(y))\right|\leq 2\varepsilon +\sup_{\substack{x_1',x_2',y_1',y_2'\in F':\\d_{F'}(x_1',x_2'),d_{F'}(y_1',y_2')\leq 7\varepsilon}}\sup_{t\in I}\left|q_t'(x_1',y_1')-q_t(x_2',y_2')\right|,\]
and so $q_t'\circ f=q_t$ for every $t\in I$ follows from the continuity properties of $q'$. This completes the proof of the fact that: if $\Delta_I\left((F,\pi,q),(F',\pi',q')\right)=0$, then the triples $(F,\pi,q)$ and $(F',\pi',q')$ are equivalent in the sense described at the start of the section. Consequently, $\Delta_I$ is indeed positive definite on the set of equivalence classes $\mathcal{M}_I$.

For the triangle inequality, we closely follow the proof of \cite[Lemma 5.2]{GPW}. Let $(F^{(i)},\pi^{(i)},q^{(i)})$ be an element of $\mathcal{M}_I$, $i=1,2,3$. Suppose that $\Delta_I((F^{(1)},\pi^{(1)},q^{(1)}),(F^{(2)},\pi^{(2)},q^{(2)}))< \delta_1$, so that we can find a metric space $Z_1$, isometric embeddings $\phi_{1,1}:F^{(1)}\rightarrow Z_1$ and $\phi_{2,1}:F^{(2)}\rightarrow Z_1$ and correspondence $\mathcal{C}_1$ between $F^{(1)}$ and $F^{(2)}$ such that the sum of quantities in the defining infimum of $\Delta_I$ is bounded above by $\delta_1$. If
$\Delta_I((F^{(2)},\pi^{(2)},q^{(2)}),(F^{(3)},\pi^{(3)},q^{(3)}))<\delta_2$, we define $Z_2$ ,$\phi_{2,2}$, $\phi_{3,2}$, $\mathcal{C}_2$ in an analogous way. Now, set $Z$ to be the disjoint union of $Z_1$ and $Z_2$, and define a distance on it by setting $d_Z|_{Z_i\times Z_i}=d_{Z_i}$ for $i=1,2$, and for $x\in Z_1$, $y\in Z_2$,
\[d_{Z}(x,y):=\inf_{z\in F^{(2)}}\left(d_{Z_1}(x,\phi_{2,1}(z))+d_{Z_2}(\phi_{2,2}(z),y)\right).\]
Abusing notation slightly, it is then the case that, after points separated by a 0 distance have been identified, $(Z,d_Z)$ is a metric space into which there is a natural isometric embedding $\phi_i$ of $Z_i$, $i=1,2$. In this space, we have that
\begin{eqnarray*}
\lefteqn{
d_H^Z(\phi_1(\phi_{1,1}(F^{(1)})),\phi_2(\phi_{3,2}(F^{(3)})))}\\
&\leq& d_H^{Z_1}(\phi_{1,1}(F^{(1)}),\phi_{2,1}(F^{(2)}))+d_H^{Z_2}(\phi_{2,2}(F^{(2)}),\phi_{3,2}(F^{(3)})),
\end{eqnarray*}
where we have applied the fact that $\phi_1(\phi_{2,1}(y))=\phi_2(\phi_{2,2}(y))$ for every $y\in F^{(2)}$, and so $\phi_1(\phi_{2,1}(F^{(2)}))=\phi_2(\phi_{2,2}(F^{(2)}))$ as subsets of $Z$. A similar bound applies to the embedded measures. Now, let
\[\mathcal{C}:=\{(x,z)\in F^{(1)}\times F^{(3)}:\exists y\in F^{(2)}\mbox{ such that } (x,y)\in\mathcal{C}_1, (y,z)\in \mathcal{C}_2\},\]
then if $(x,z)\in \mathcal{C}$,
\[d_Z(\phi_1(\phi_{1,1}(x)),\phi_2(\phi_{3,2}(z)))\leq d_{Z_1}(\phi_{1,1}(x),\phi_{2,1}(y))+d_{Z_2}(\phi_{2,2}(y),\phi_{3,2}(z)),\]
where $y\in F^{(2)}$ is chosen such that $(x,y)\in\mathcal{C}_1$ and $(y,z)\in \mathcal{C}_2$, and we again note $\phi_1(\phi_{2,1}(y))=\phi_2(\phi_{2,2}(y))$. Proceeding in the same fashion, one can deduce a corresponding bound involving $q^{(i)}$, $i=1,2,3$. Putting these pieces together, it is elementary to deduce that
\[\Delta_I((F^{(1)},\pi^{(1)},q^{(1)}),(F^{(3)},\pi^{(3)},q^{(3)}))\leq \delta_1+\delta_2,\]
and the triangle inequality follows. Thus we have proved that $(\mathcal{M}_I,\Delta_I)$ is a metric space.

To complete the proof, we only need to show separability. This is straightforward, however, as for any element of $\mathcal{M}_I$, one can construct an approximating sequence that incorporates only: metric spaces with a finite number of points and rational distances between them, probability measures on these with a rational mass at each point, and functions that are defined (at each coordinate pair) to be equal to rational values at a finite collection of rational time points and are linear between these. To be more explicit, let $(F,\pi,q)$ be an element of $\mathcal{M}_I$, and then define a sequence $(F^N,\pi^N,q^N)_{N\geq 1}$ as follows. First, let $F^N$ be a finite $N^{-1}$-net of $F$, which exists because $F$ is compact. By perturbing $d_F$, it is possible to define a metric $d_{F^N}$ on $F^N$ such that $|d_{F^N}(x,y)-d_F(x,y)|\leq N^{-1}$ and moreover $d_{F^N}(x,y)\in \mathbb{Q}$ for all $x,y\in F^N$. Now, since $F^N$ is an $N^{-1}$-net of $F$, it is possible to choose a partition $(A_x)_{x\in F^N}$ of $F$ such that $x\in A_x$ and the diameter of $A_x$ (with respect to $d_F$) is no greater than $2N^{-1}$. Moreover, it is possible to choose the partition in such a way that $A_x$ is measurable for each $x\in F^N$. We construct a probability measure on $F^N$ by choosing $\pi^N(\{x\})\in \mathbb{Q}$ such that $|\pi^N(\{x\})-\pi(A_x)|\leq N^{-1}$ (subject to the constraint that $\sum_{x\in F_N}\pi^N(\{x\})=1)$. Finally, define $\varepsilon_N$ by setting
\[\varepsilon_N:=\sup_{\substack{s,t\in I:\\|s-t|\leq N^{-1}}}\sup_{\substack{x,x',y,y'\in F:\\d_F(x,x'),d_F(y,y')\leq N^{-1}}}\left|q_s(x,y)-q_t(x',y')\right|,\]
so that, by the joint continuity of $q$, $\varepsilon_N\rightarrow 0$ as $N\rightarrow\infty$. Let $\inf I\leq t_0\leq t_1\leq \dots\leq t_K\leq \sup I$ be a set of rational times such that $|t_0-\inf I|$, $|\sup I- t_K|$, $|t_{i+1}-t_i|\leq N^{-1}$, choose $q^N_{t_i}(x,y)\in \mathbb{Q}$ such that $|q^N_{t_i}(x,y)-q_{t_i}(x,y)|\leq N^{-1}$ for each $x,y\in F^N$, and then extend $q^N$ to have domain $F^N\times F^N\times I$ by linear interpolation in $t$ at each pair of vertices. This construction readily yields that $\Delta_I((F,\pi,q),(F^N,\pi^N,q^N))\leq 6N^{-1}+3\varepsilon_N \rightarrow 0$. Since the class of triples from which the approximating sequence is chosen is clearly countable, this completes the proof of separability.
\end{proof}

We will say that a sequence in $\mathcal{M}_I$ converges in a spectral Gromov-Hausdorff sense if it converges to a limit in this space with respect to the metric $\Delta_I$. We note that in the framework of compact Riemannian manifolds, different but related notions of spectral distances were introduced
by B\'erard, Besson and Gallot (\cite{BBG}) and by Kasue and Kumura (\cite{KasK}). Moreover, by applying our characterization of spectral Gromov-Hausdorff convergence, we are able to deduce that if Assumption \ref{assu1} holds, then we can isometrically embed all the rescaled graphs, measures and transition densities upon them into a common metric space $(E,d_E)$ so that they converge to the relevant limit objects in a more standard way, as the following lemma makes precise. Note that in the proof of the result and henceforth we define balls in the space $(E,d_E)$ by setting $B_E(x,r):=\{x\in E:d_E(x,y)<r\}$.

{\lem\label{embed} Suppose that Assumption \ref{assu1} is satisfied. For any compact interval $I\subset(0,\infty)$, there exist isometric embeddings of $(V(G^N),d_{G^N})$, $N\geq1$, and $(F,d_F)$ into a common metric space $(E,d_E)$ such that
\begin{equation}\label{a}
\lim_{N\rightarrow\infty}
d^E_H(V(G^N),F)= 0,
\end{equation}
\begin{equation}\label{b}
\lim_{N\rightarrow\infty}
d_P^E(\pi^{N},\pi)=0,
\end{equation}
and also,
\begin{equation}\label{c}
\lim_{N\rightarrow\infty}\sup_{x,y\in F}\sup_{t\in I} \left|q^N_{\gamma(N)t}(g_N(x),g_N(y))-q_t(x,y)\right|=0,
\end{equation}
where, for brevity, we have identified the spaces $(V(G^N),d_{G^N})$, $N\geq1$, and $(F,d_F)$, and the measures upon them with their isometric embeddings in $(E,d_E)$. For each $x\in F$, we define $g_N(x)$ to be a vertex in $V(G^N)$ minimizing $d_{E}(x,y)$ over $y\in V(G^N)$.}
\begin{proof} Fix a compact interval $I\subset(0,\infty)$. By Assumption \ref{assu1}, for each $N\geq 1$ it is possible to find metric spaces $(E_N,d_{N})$, isometric embeddings $\phi_N:(V(G^N),d_{G^N})\rightarrow (E_N,d_{N})$, $\phi_N':(F,d_F)\rightarrow (E_N,d_{N})$ and correspondences $\mathcal{C}_N$ between $V(G^N)$ and $F$ such that, identifying the original objects and their embeddings,
\begin{eqnarray}
\lefteqn{d_H^{E^N}(V(G^N),F)+d_P^{E^N}(\pi^N,\pi)}\hspace{350pt}\nonumber\\
+\sup_{(x,x'),(y,y')\in\mathcal{C}_N}
\left(d_{N}(x,x')+d_{N}(y,y')+\sup_{t\in I}\left|q_{\gamma(N)t}^N(x,y)-q_t(x',y')\right|\right)&\leq&\varepsilon_N,\label{decay}
\end{eqnarray}
where $\varepsilon_N\rightarrow 0$. Now, proceeding similarly to the proof of the triangle inequality in Lemma \ref{ghlem}, set $E$ to be the disjoint union of $E^N$, $N\geq1$, and define a distance on it by setting $d_E|_{E^N\times E^N}=d_{N}$ for $N\geq1$, and for $x\in E^N$, $x'\in E^{N'}$, $N\neq N'$, set
\[d_{E}(x,x'):=\inf_{y\in F}\left(d_{N}(x,y)+d_{N'}(y,x')\right).\]
Quotienting out points that are separated by distance 0 results in a metric space $(E,d_E)$ (again, this is a slight abuse of notation), into which we have natural isometric embeddings of the metric spaces $(V(G^N),d_{G^N})$, $N\geq1$, and $(F,d_F)$. Moreover, in the metric space $(E,d_E)$, it readily follows from (\ref{decay}) that the relevant isometrically embedded objects satisfy (\ref{a}) and (\ref{b}). To prove (\ref{c}), first note that for every $x\in V(G^N)$, $N\geq1$, there exists an $x'\in F$ such that $(x,x')\in \mathcal{C}_N$. This implies that $d_E(x,x')\leq\varepsilon_N$, and so, for any $\delta>0$,
\begin{eqnarray}
\lefteqn{\sup_{\substack{x,y,z\in V(G^N):\\d_{G^N}(y,z)\leq \delta}}\sup_{t\in I} \left|q^N_{\gamma(N)t}(x,y)-q^N_{\gamma(N)t}(x,z)\right|}\nonumber\\
&\leq&
2\varepsilon_N+\sup_{\substack{x,y,z\in F:\\d_{F}(y,z)\leq \delta+2\varepsilon_N}}\sup_{t\in I} \left|q_{t}(x,y)-q_{t}(x,z)\right|.\label{tightbound}
\end{eqnarray}
Now, for every $x\in F$ and $N\geq 1$, there exists an $x'\in V(G^N)$ such that $(x',x)\in \mathcal{C}_N$, and so $d_E(x',x)\leq\varepsilon_N$. Therefore, since $g_N(x)$ is the closest vertex of $V(G^N)$ to $x$,
\[g_N(x)\in B_E(x,2\varepsilon_N)\cap V(G^N)\subseteq B_E(x',3\varepsilon_N)\cap V(G^N)=B_{V(G^N)}(x',3\varepsilon_N).\]
Consequently,
\begin{eqnarray*}
\lefteqn{\sup_{x,y\in F}\sup_{t\in I} \left|q^N_{\gamma(N)t}(g_N(x),g_N(y))-q_t(x,y)\right|}\\
&\leq&
\varepsilon_N+2\sup_{\substack{x,y,z\in V(G^N):\\d_{G^N}(y,z)\leq 3\varepsilon_N}}\sup_{t \in I}\left|q_{\gamma(N)t}^{N}(x,y)-q_{\gamma(N)t}^{N}(x,z)\right|\\
&\leq &5\varepsilon_N+2\sup_{\substack{x,y,z\in F:\\d_{F}(y,z)\leq 5\varepsilon_N}}\sup_{t\in I} \left|q_{t}(x,y)-q_{t}(x,z)\right|,
\end{eqnarray*}
where the second inequality is an application of (\ref{tightbound}). Letting $N\rightarrow\infty$ and applying the joint continuity of $(q_t(x,y))_{x,y\in F,t>0}$, we obtain the desired result.
\end{proof}

For our later convenience, let us note a useful tightness condition for the rescaled transition densities that was essentially established in the proof of the previous result.

{\lem\label{assu2} Suppose that Assumption \ref{assu1} holds. For any compact interval $I\subset (0,\infty)$,
\begin{equation}\label{tightcond}
\lim_{\delta\rightarrow 0}\limsup_{N\rightarrow\infty}\sup_{\substack{x,y,z\in V(G^N):\\d_{G^N}(y,z)\leq \delta}}\sup_{t \in I}\left|q_{\gamma(N)t}^{N}(x,y)-q_{\gamma(N)t}^{N}(x,z)\right|=0.
\end{equation}}
\begin{proof} Recalling the continuity property of $q$, taking the limit as $N\rightarrow\infty$ in (\ref{tightbound}) yields
\[\limsup_{N\rightarrow\infty}\sup_{\substack{x,y,z\in V(G^N):\\d_{G^N}(y,z)\leq \delta}}\sup_{t \in I}\left|q_{\gamma(N)t}^{N}(x,y)-q_{\gamma(N)t}^{N}(x,z)\right|\leq \sup_{\substack{x,y,z\in F:\\d_{F}(y,z)\leq \delta}}\sup_{t \in I}\left|q_t(x,y)-q_t(x,z)\right|.\]
Again appealing to the continuity of $q$, the right-hand side here converges to 0 as $\delta\rightarrow 0$, which completes the proof.
\end{proof}

It is straightforward to reverse the conclusions of the previous two lemmas to check that if (\ref{a}), (\ref{b}), (\ref{c}) and (\ref{tightcond}) hold, then so does Assumption \ref{assu1}. Indeed, under these assumptions, we have isometric embeddings of $(V(G^N),d_{G^N})$, $N\geq 1$, and $(F,d_F)$ into a common metric space $(E,d_E)$ for which: (\ref{a}) gives the Hausdorff convergence of sets; (\ref{b}) gives the Prohorov convergence of measures; and moreover, it is elementary to check from (\ref{c}) and (\ref{tightcond}) that, with respect to the correspondences
\[\mathcal{C}_N:=\left\{(x,x')\in F\times V(G^N):d_E(x,x')\leq N^{-1}\right\},\]
the relevant transition densities converge uniformly, as described in the definition of the metric $\Delta_I$.  Thus, in examples, it will suffice to check these equivalent conditions when seeking to verify Assumption \ref{assu1}. In fact, it is further possible to weaken these assumptions slightly by appealing to a local limit theorem from \cite{CHLLT}. To be precise, because we are assuming that the transition densities of the graph satisfy the tightness condition of (\ref{tightcond}), we can apply \cite[Theorem 15]{CHLLT}, to replace the local convergence statement of (\ref{c}) with a central limit-type convergence statement. Note that, although in \cite{CHLLT} it was assumed that the metric on $G^N$ was a shortest path graph distance, exactly the same argument yields the corresponding conclusion in our setting, and so we simply state the result.

{\propn [cf. {\cite[Theorem 15]{CHLLT}}]\label{local} Suppose that $(V(G^N), d_{G^N})$, $N\geq 1,$ and $(F,d_F)$ can be isometrically embedded into a common metric space $(E,d_E)$ in such a way that (\ref{a}) and (\ref{b}) are both satisfied. Moreover, assume that there exists a dense subset $F^*$of $F$ such that, for any compact interval $I\subset (0,\infty)$, $x\in F^*$, $y\in F$, $r>0$,
\begin{equation}\label{CLT-dist}
\lim_{N\rightarrow\infty}\mathbf{P}^{G^N}_{g_N(x)}\left(X^{G^N}_{\lfloor \gamma(N)t\rfloor}
\in B_E(y,r)\right)=
\mathbf{P}^F_{x}\left(X^{F}_t\in B_E(y,r)\right)\end{equation}
uniformly for $t\in I$, and also (\ref{tightcond}) holds. Then Assumption 1 holds.}
\bigskip

To complete this section, let us observe that \cite{CHLLT} also provides two ways to check (\ref{tightcond}): one involving a resistance estimate on the graphs in the sequence (\cite[Proposition 17]{CHLLT}), and one involving the parabolic Harnack inequality (\cite[Proposition 16]{CHLLT}). Since the first of these two methods will be applied in several of our examples later, let us recall the result here. To allow us to state the result, we define $R_{G^N}(x,y)$ to be the
resistance between $x$ and $y$ in $V(G^N)$ (see \eqref{eq:residef}),
when we suppose that $G^N$ is an electrical network with conductances of edges being given by the weight function $\mu^{G^N}$. This defines a metric on $V(G^N)$, for which the following result is proved as \cite[Proposition 17]{CHLLT}. As above, note that although it was a shortest path graph distance considered in \cite{CHLLT}, the same proof applies for a general distance on the graph in question. Moreover, the statement of the lemma is slightly different from that of the corresponding result in \cite{CHLLT}, because there the scaling $\alpha(n)$ was absorbed into the definition of the metric.

{\lem[cf. {\cite[Proposition 17]{CHLLT}}]\label{tightlem} Suppose that there exists a sequence $(\alpha(N))_{N\geq 1}$ and constants $\kappa, c_1,c_2, c_3\in (0,\infty)$ such that
\[R_{G^N}(x,y)\leq c_1\left(\alpha(N)d_{G^N}(x,y)\right)^{\kappa},\hspace{20pt}\forall x,y\in V(G^N),\]
and also
\[c_2\gamma(N)\leq \alpha(N)^\kappa\beta(N)\leq c_3\gamma(N),\]
where $\beta(N):=\sum_{x,y\in V(G^N)}\mu^{G^N}_{xy}$, then (\ref{tightcond}) holds.}

\section{Convergence of $L^p$-mixing times}\label{mixingtime}

\subsection{Proof of Theorem \ref{main}}

In this subsection we prove the mixing time convergence result of Theorem \ref{main}. Throughout, we will suppose that Assumption \ref{assu1} holds and that the graphs $G^N$ and limiting metric space $F$ have been embedded into a common metric space $(E,d_E)$ in the way described by Lemma \ref{embed}.

Recall from the introduction the definition of $D_p(x,t)=\|q_t(x,\cdot)-1\|_{L^p(\pi)}$, the $L^p$-distance from stationarity of the process $X^F$ started from $x$ at time $t$. By applying the continuity of $(q_t(x,y))_{x,y\in F,t>0}$, compactness of $F$ and finiteness of $\pi$, it is easy to check that this quantity is finite for every $x\in F$ and $t>0$. The next lemma collects together a number of other basic properties of $D_p(x,t)$ that we will apply later (the first part is a minor extension of \cite[Proposition 3.1]{S-C}, in our setting).

{\lem\label{dplem} Let $p\in[1,\infty]$. For every $x\in F$, the function $t\mapsto D_p(x,t)$ is continuous and strictly
decreasing. Furthermore, we have
\begin{equation}\label{0lim}
\lim_{t\rightarrow 0}D_p(x,t)\geq 2.
\end{equation}}
\begin{proof} That the function $t\mapsto D_p(x,t)$ is continuous is clear from
\eqref{asum-b}, and so we turn to checking that it is strictly decreasing. First, a standard argument involving an application
of Jensen's inequality and the invariance of $\pi$ allows one to deduce that $\|P_t f\|_{L^p(\pi)}\leq \|f\|_{L^p(\pi)}$
for any $f\in L^p(F,\pi)$, where $(P_t)_{t\geq0}$ is the semigroup naturally associated with the transition density
$(q_t(x,y))_{x,y\in F,t>0}$.  Now, suppose $f\in L^p(F,\pi)$ is such that $\|P_t f\|_{L^p(\pi)}=\|f\|_{L^p(\pi)}$,
and define $f_1(y):=|P_tf(y)|^p$ and $f_2(y):=P_t(|f|^p)(y)$. By the assumption on $f$ and the fact that $X^F$ is conservative and $\pi$-symmetric, we have that
\begin{eqnarray*}
\int_Ff_1 d\pi&=&\int_F |P_tf(y)|^p\pi(dy)\\
&=&\int_F |f(y)|^p\pi(dy)\\
&=&\int_F |f(y)|^p \int_Fq_t(y,z)\pi(dz)\pi(dy)\\
&=&\int_F \int_F |f(y)|^p q_t(z,y)\pi(dy)\pi(dz)\\
&=&\int_F P_t(|f|^p)(z)\pi(dz)\\
&=&\int_Ff_2d\pi.
\end{eqnarray*}
Furthermore, Jensen's inequality implies $f_1(y)\leq f_2(y)$. Thus, it must be the
 case that $f_1(y)=f_2(y)$, $\pi$-a.e. In particular, because $\pi$ is a probability measure, there exists a $y\in F$
 such that $f_1(y)=f_2(y)$.

 In the case $p>1$, the conclusion of the previous paragraph readily implies that $f$ is constant $q_t(y,z)\pi(dz)$-a.e.
 Recalling the assumption that $q_t(y,z)>0$ everywhere, namely \eqref{asum-c}, it must therefore hold that
 $f$ is constant $\pi$-a.e. Observing that for $s,t>0$ we can write
$D_p(x,s+t)=\|P_s(q_t(x,\cdot)-1)\|_{L^p(\pi)}$, it follows that $D_p(x,s+t)<D_p(x,t)$ if and only if $q_t(x,\cdot)=1$, $\pi$-a.e. However, condition \eqref{asum-d} and the assumption that the transition density is continuous imply that there exists a non-empty open set on which $q_t(x,\cdot)\neq 1$. Thus, because $\pi$ has full support, it is not the case that $q_t(x,\cdot)=1$, $\pi$-a.e., and we must have $D_p(x,s+t)<D_p(x,t)$, as desired.

For $p=1$, the result
$f_1(y)=f_2(y)$ implies that $f$ is either non-negative or non-positive, $\pi$-a.e. Consequently, if we suppose that $D_p(x,s+t)=\|P_s(q_t(x,\cdot)-1)\|_{L^p(\pi)}=D_p(x,t)$ for some $s>0$, then it must be the case that $q_t(x,\cdot)-1$ is either non-negative or non-positive. However, since $\int_F(q_t(x,y)-1)\pi(dy)=0$ (due to \eqref{asum-a}) and \eqref{asum-d} holds, we arrive at a contradiction. In particular, it must be the case that $D_p(x,s+t)<D_p(x,t)$, and this completes the proof of strict monotonicity.

To establish the limit in (\ref{0lim}), it will suffice to prove the result in the case $p=1$ (obtaining the result
for other values of $p$ is then simply Jensen's inequality). Let $x\in F$ and $r>0$, then
\begin{eqnarray*}
D_1(x,t)&\geq&\int_{B_E(x,r)}(q_t(x,y)-1)\pi(dy)+\int_{B_E(x,r)^c}(1-q_t(x,y))\pi(dy)\\
&=&2\mathbf{P}_x\left(X^F_t\in B_E(x,r)\right)-2\pi(B_E(x,r)),
\end{eqnarray*}
where \eqref{asum-a} is used in the last equality. Since $X^F$ is a Hunt process, the first term here converges to 2
as $t\rightarrow 0$. Furthermore, because $\pi$ is non-atomic, the second term can be made arbitrarily small by
suitable choice of $r$. The result follows.
\end{proof}

We continue by defining the $L^p$-mixing time at $x\in F$ by setting
\[t_{\rm mix}^p(x):=\inf\{t>0:D_p(x,t)\leq 1/4\}.\]
In fact, the previous lemma yields that $t_{\rm mix}^p(x)$ is the unique value of $t\in(0,\infty)$ such that $D_p(x,t)= 1/4$ (when (\ref{erg}) holds at $x$). Similarly, define the $L^p$-mixing time of $x\in V(G^N)$ by setting
\[t_{\rm mix}^{N,p}(x):=\inf\{t>0:D^{N}_p(x,t)\leq 1/4\},\]
where $D^{N}_p(x,m)=\|q^N_m(x,\cdot)-1\|_{L^p(\pi^N)}$. That the discrete mixing times at a point converge when suitably rescaled to the continuous mixing time there is the conclusion of the following proposition.

{\propn \label{point} Suppose that Assumption \ref{assu1} is satisfied. If $p\in[1,\infty]$ is such that (\ref{erg}) holds for $x\in F$, then
\[\lim_{N\rightarrow\infty}\gamma(N)^{-1}t_{\rm mix}^{N,p}(g_N(x))=t_{\rm mix}^p(x),\]
where, as in the statement of Lemma \ref{embed}, $g_N(x)$ is a vertex in $V(G^N)$ that minimizes the distance $d_E(x,y)$ over $V(G^N)$.}
\begin{proof} Suppose $p\in[1,\infty]$ is such that (\ref{erg}) holds for $x\in F$, set $t_0:=t_{\rm mix}^p(x)\in (0,\infty)$, and fix $\varepsilon>0$. By
\eqref{asum-b} and the tightness of Lemma \ref{assu2}, there exists a $\delta>0$ such that
\begin{equation}\label{cont}
\sup_{t\in I}\sup_{\substack{y,z\in F:\\d_E(y,z)<2\delta}}\left||q_t(x,y)-1|^p-|q_t(x,z)-1|^p\right|<\varepsilon,
\end{equation}
\begin{equation}\label{discont}
\limsup_{N\rightarrow\infty}\sup_{t\in I}\sup_{\substack{y,z\in V(G^N):\\d_{G^N}(y,z)<3\delta}}\left||q^N_{\gamma(N)t}(g_N(x),y)-1|^p-|q^N_{\gamma(N)t}(g_N(x),z)-1|^p\right|<\varepsilon,
\end{equation}
where $I:=[t_0/2,2t_0]$. Moreover, by the compactness of $F$, there exists a finite collection of balls $(B_E(x_i,\delta))_{i=1}^k$ covering $F$. Define $A_1:=B(x_1,2\delta)$, and $A_i:=B_E(x_i,2\delta)\backslash\cup_{j=1}^{i-1}B_E(x_i,2\delta)$ for $i=2,\dots,k$, so that $(A_i)_{i=1}^k$ is a disjoint cover of the $\delta$-enlargement of $F$.

We observe
\[|D_p(x,t)^p-D_p^N(g_N(x),\gamma(N)t)^p|\leq T_1+T_2+T_3+T_4,\]
where
\begin{eqnarray*}
T_1&:=& \left|\int_F|q_t(x,y)-1|^p\pi(dy)-\sum_{i=1}^k|q_t(x,x_i)-1|^p\pi(A_i)\right|,\\
T_2&:=& \left|\sum_{i=1}^k|q_t(x,x_i)-1|^p\pi(A_i)-\sum_{i=1}^k|q_t(x,x_i)-1|^p\pi^N(A_i)\right|,\\
T_3&:=& \left|\sum_{i=1}^k|q_t(x,x_i)-1|^p\pi^N(A_i)-\sum_{i=1}^k|q^N_{\gamma(N)t}(g_N(x),g_N(x_i))-1|^p\pi^N(A_i)\right|,\\
T_4&:=& \left|\sum_{i=1}^k|q^N_{\gamma(N)t}(g_N(x),g_N(x_i))-1|^p\pi^N(A_i)-\int_{V(G^N)}|q^N_{\gamma(N)t}(g_N(x),y)-1|^p\pi^N(dy)\right|.
\end{eqnarray*}
Now, suppose $t\in I$. From (\ref{cont}), we immediately deduce that $T_1\leq\varepsilon$. For $T_2$, we first observe that the fact balls are $\pi$-continuity sets implies that $A_1,\dots, A_k$ are also $\pi$-continuity sets. Hence $\pi^N(A_i)\rightarrow \pi(A_i)$ for each $i=1,\dots,k$, and so $T_2\leq \varepsilon$ for large $N$. That $T_3\leq \varepsilon$ for large $N$ is a straightforward consequence of Lemma \ref{embed}. Finally, applying the fact that $d^E_H(F,V(G^N))\rightarrow 0$, we deduce that, for large $N$, $(A_i)_{i=1}^k$ is a disjoint cover for $V(G^N)$. Since $g_N(x_i)\in B_E(x_i,\delta)$ for large $N$, we also have that $d_{G^N}(y,g_N(x_i))\leq 3\delta$, uniformly over $y\in A_i$, $i=1,\dots,k$. Thus we can appeal to (\ref{discont}) to deduce that it is also the case that $T_4\leq \varepsilon$ for large $N$. In fact, each of these bounds can be assumed to hold uniformly over $t\in I$, thereby demonstrating that
\begin{equation}\label{dconv}
\lim_{N\rightarrow\infty}\sup_{t\in I}\left|D_p(x,t)-D_p^N(g_N(x),\gamma(N)t)\right|=0.
\end{equation}
Since $t\mapsto D_p^N(g_N(x),\gamma(N)t)$ is a decreasing function in $t$ for every $N$ (cf. \cite[Proposition 3.1]{S-C}) and $t\mapsto D_p(x,t)$ is strictly decreasing, the proposition follows.
\end{proof}

{\rem {\rm In the case $p=2$, the proof of the previous result greatly simplifies. In particular, we note that
\begin{equation}\label{l2obs}
D_2(x,t)^2=\|q_t(x,\cdot)-1\|_2^2=q_{2t}(x,x)-1,
\end{equation}
and similarly
\[D_2^N(x,\gamma(N)t)^2=\|q^N_{\gamma(N)t}(x,\cdot)-1\|_2^2=q^N_{2{\gamma(N)t}}(x,x)-1.\]
Hence the limit at (\ref{dconv}) is an immediate consequence of the local limit result of (\ref{c}), and we do not have to concern ourselves with estimating the relevant integrals directly.}}
\bigskip

To extend the above proposition to the corresponding result for the mixing times of the entire spaces, we will appeal to the following lemma, which establishes a continuity property for the  $L^p$-mixing times from fixed starting points in the limiting space, and a related tightness property for the discrete approximations.

{\lem\label{lem2} Suppose $p\in[1,\infty]$ is such that (\ref{erg}) holds for $x\in F$, then the following statements are true.\\
(a) The function $y\mapsto t_{\rm mix}^p(y)$ is continuous at $x$.\\
(b) Under Assumption \ref{assu1}, it is the case that
\[\lim_{\delta\rightarrow0}\limsup_{N\rightarrow\infty}\sup_{\substack{y\in V(G^N):\\d_{G^N}(g_N(x),y)<\delta}}\gamma(N)^{-1}\left|t_{\rm mix}^{N,p}(y)-t_{\rm mix}^{N,p}(g_N(x))\right|=0.\]}
\begin{proof} Consider $p\in[1,\infty]$ such that (\ref{erg}) holds for $x\in F$, so that $t_0:=t_{\rm mix}^p(x)$ is finite, and let $\varepsilon\in(0,t_0/2)$. Since the function $t\mapsto D_p(x,t)$ is strictly decreasing (by Lemma \ref{dplem}), there exists an $\eta>0$ such that $D_p(x,t_0-\varepsilon)>D_p(x,t_0)+\eta=1/4+\eta$ and also
$D_p(x,t_0+\varepsilon)<1/4-\eta$. By the continuity of $(q_t(x,y))_{x,y\in F,t>0}$, there also exists a $\delta>0$ such that
\[\sup_{t\in [t_0-\varepsilon,t_0+\varepsilon]}\sup_{\substack{y\in F:\\d_{F}(x,y)<\delta}}\left|D_p(x,t)-D_p(y,t)\right|<\eta.\]
Hence if $y\in B_F(x,\delta)$, then
\[D_p(y,t_0-\varepsilon)>D_p(x,t_0-\varepsilon)-\eta>\frac{1}{4},\]
\[D_p(y,t_0+\varepsilon)<D_p(x,t_0+\varepsilon)+\eta<\frac{1}{4}.\]
This implies that $t_{\rm mix}^{p}(y)\in[t_0-\varepsilon,t_0+\varepsilon]$, and (a) follows.

The proof of part (b) is similar. In particular, choose $\eta$ as above and note that (\ref{dconv}) implies that $D_p^N(g_N(x),\gamma(N)(t_0-\varepsilon))>1/4+\eta/2$ and $D_p^N(g_N(x),\gamma(N)(t_0+\varepsilon))<1/4-\eta/2$ for large $N$. Furthermore, by the transition density tightness of Lemma \ref{assu2}, there exists a $\delta>0$ such that
\[\sup_{t\in [t_0-\varepsilon,t_0+\varepsilon]}\sup_{\substack{y\in V(G^N):\\d_{G^N}(g_N(x),y)<\delta}}\left|D_p^N(g_N(x),\gamma(N)t)-D_p^N(y,\gamma(N)t)\right|<\frac{\eta}{2},\]
for large $N$. Hence if $N$ is large and $y\in V(G^N)$ is such that $d_{G^N}(g_N(x),y)<\delta$, then $D_p^N(y,\gamma(N)(t_0-\varepsilon))>1/4$, and $D_p^N(y,\gamma(N)(t_0+\varepsilon))<1/4$. This implies that $\gamma(N)^{-1}t_{\rm mix}^{N,p}(y)\in[t_0-\varepsilon,t_0+\varepsilon]$. Since it is trivially true that, once $N$ is large enough, this result can be applied with $y=g_N(x)$, the result follows.
\end{proof}

We are now ready to give the proof of our main result.

\begin{proof}[Proof of Theorem \ref{main}] Observe that, under the assumptions of the theorem, Lemma \ref{lem2}(a) implies that the function $(t_{\rm mix}^p(x))_{x\in F}$ is continuous. Since $F$ is compact, the supremum of  $(t_{\rm mix}^p(x))_{x\in F}$ is therefore finite. Now, it is an elementary exercise to check that we can write the $L^p$-mixing time of $F$, as defined at (\ref{mixingdef}), in the following way:
\begin{equation}\label{supexp}
t_{\rm mix}^p(F)=\sup_{x\in F}t_{\rm mix}^p(x).
\end{equation}
Consequently $t_{\rm mix}^p(F)\in (0,\infty)$, as desired.

To complete the proof, we are required to demonstrate the convergence statement of (\ref{mainconv}). Fix $\varepsilon>0$. For every $x\in F$, Proposition \ref{point} and Lemma \ref{lem2}(b) allow us to choose $\delta(x)>0$ and $N(x)<\infty$ such that
\[\sup_{N\geq N(x)}\left|\gamma(N)^{-1}t_{\rm mix}^{N,p}(g_N(x))-t_{\rm mix}^p(x)\right|\leq \varepsilon,\]
\[\sup_{N\geq N(x)}\sup_{\substack{y\in V(G^N):\\d_{G^N}(g_N(x),y)<4\delta(x)}}\gamma(N)^{-1}\left|t_{\rm mix}^{N,p}(g_N(x))-t_{\rm mix}^{N,p}(y)\right|\leq\varepsilon.\]
Since $(B_E(x,\delta(x)))_{x\in F}$ is an open cover for $F$, by compactness it admits a finite subcover $(B_E(x,\delta(x)))_{x\in \mathcal{X}}$. Moreover, because $d^E_H(F,V(G^N))\rightarrow 0$, there exists an $N_0>0$ such that if $N\geq N_0$, then $(B_E(x,2\delta(x)))_{x\in \mathcal{X}}$ is a cover for $V(G^N)$. Applying this choice of $\mathcal{X}$, we have for $N\geq N_0\vee\max_{x\in\mathcal{X}}N(x)$ that
\[\gamma(N)^{-1}t_{\rm mix}^{p}(G^N)\leq\sup_{x\in \mathcal{X}}\gamma(N)^{-1}t_{\rm mix}^{N,p}(g_N(x))+\varepsilon
\leq \sup_{x\in \mathcal{X}}t_{\rm mix}^{p}(x)+2\varepsilon\leq  t_{\rm mix}^{p}(F)+2\varepsilon,\]
where we note that, similarly to (\ref{supexp}), the $L^p$-mixing time of the graph $G^N$ can be written as
\[t_{\rm mix}^{p}(G^N)=\sup_{x\in V(G^{N})}t_{\rm mix}^{N,p}(x).\]
Furthermore, if $x_0\in F$ is chosen such that $t_{\rm mix}^{p}(x_0)\geq t_{\rm mix}^{p}(F)-\varepsilon$, then, for large $N$,
\[\gamma(N)^{-1}t_{\rm mix}^{N,p}(G^N)\geq \gamma(N)^{-1}t_{\rm mix}^{N,p}(g_N(x_0))\geq t_{\rm mix}^{p}(x_0)-\varepsilon\geq   t_{\rm mix}^{p}(F)-2\varepsilon,\]
where we have again made use of Proposition \ref{point}. Since $\varepsilon>0$ was arbitrary, we are done.
\end{proof}

\subsection{Distinguished starting points}\label{fixed}

In certain situations, convergence of transition densities might only be known with respect to a single distinguished starting point. This is the case, for instance, in two of the most important examples we present in Section \ref{examsec} -- critical Galton-Watson trees and the critical Erd\H{o}s-R\'{e}nyi random graph. In such settings, it is only possible to prove a convergence result for the mixing time from the distinguished point. It is the purpose of this subsection to present a precise conclusion of this kind.

Consider, for a compact interval $I\subset(0,\infty)$, the space of triples of the form $(F,\pi,q)$, where $F=(F,d_F,\rho)$ is a non-empty compact metric space with distinguished vertex $\rho$, $\pi$ is a Borel probability measure on $F$ and $q=(q_t(x,y))_{x,y\in F, t\in I}$ is a jointly continuous real-valued function of $(t,x,y)$; this is the same as the collection $\tilde{\mathcal{M}}_{I}$ defined in Section \ref{gh}, though we have added the supposition that the metric spaces are pointed. We say two such elements, $(F,\pi,q)$ and $(F',\pi',q')$, are equivalent if there exists an isometry $f:F\rightarrow {F'}$ such that $f(\rho)=\rho'$, $\pi\circ f^{-1}=\pi'$ and $q_t'\circ f=q_t$ for every $t\in I$. By following the proof of Lemma \ref{ghlem}, one can check that it is possible to define a metric on the equivalence classes of this relation by simply including in the definition of $\Delta_I$ the condition that the correspondence $\mathcal{C}$ must contain $(\rho,\rho')$. We define convergence in a spectral pointed Gromov-Hausdorff sense to be with respect to this metric. The distinguished starting point version of Assumption \ref{assu1} is then as follows.

{\assu\label{assu3} Let $(G^N)_{N\geq 1}$ be a sequence of finite connected graphs with at least two vertices and one, $\rho^N$ say, distinguished, for which there exists a sequence $(\gamma(N))_{N\geq 1}$ such that, for any compact interval $I\subset (0,\infty)$,
\[{\left(\left(V(G^N),d_{G^N},\rho^N\right), \pi^N, \left(q^N_{\gamma(N)t}(\rho^N,x)\right)_{x\in V(G^N),t\in I}\right)}\]
converges to $(\left(F,d_F,\rho\right),\pi,(q_t(\rho,x))_{x\in F,t\in I})$ in a spectral pointed Gromov-Hausdorff sense, where $\rho$ is a distinguished point in $F$.}
\bigskip

The following result can then be proved in an almost identical fashion to Proposition \ref{point}, simply replacing $g_N(x)$ by $\rho^N$ and $x$ by $\rho$. In doing this, it is useful to note that if Assumption \ref{assu1} is replaced by Assumption \ref{assu3}, then we are able to include in the conclusions of Lemma \ref{embed} that $\rho^N$ converges to $\rho$ in $E$.

{\thm \label{pointthm}Suppose that Assumption \ref{assu3} is satisfied. If $p\in[1,\infty]$ is such that
(\ref{erg}) holds for $x=\rho$, then
\[\gamma(N)^{-1}t_{\rm mix}^{N,p}(\rho^N)\rightarrow t_{\rm mix}^{p}(\rho).\]}

\section{Convergence to stationarity of the transition density}

Before continuing to present example applications of the mixing time convergence results proved so far,
we describe how to check the $L^p$ convergence to stationarity of the transition density of $X^F$ in the
case when we have a spectral decomposition for it and a spectral gap. In the same setting, we will also
explain how to check the non-triviality conditions on the transition density that were made in the introduction.

Write the generator of the conservative Hunt process $X^F$ as $-\Delta$, and suppose that $\Delta$ has
a compact resolvent. Then there exists a complete orthonormal basis of $L^2(F,\pi)$, $(\varphi_k)_{k\geq 1}$
say, such that $\Delta\varphi_k=\lambda_k\varphi_k$ for all $k\geq 0$, $0\leq \lambda_0\leq \lambda_1\leq\dots$
and $\lim_{k\rightarrow\infty}\lambda_k=\infty$. By expanding as a Fourier series, we can consequently
write the transition density of $X^F$ as
\begin{eqnarray*}
q_t(x,y)&=&\sum_{k\geq 0}\left(\int_{F}q_t(x,z)\varphi_k(z)\pi(dz)\right)\varphi_k(y)\\
&=&\sum_{k\geq 0}P^F_t\varphi_k(x)\varphi_k(y)\\
&=&\sum_{k\geq 0}e^{-\lambda_kt}\varphi_k(x)\varphi_k(y),
\end{eqnarray*}
where $(P^F_t)_{t\geq0}$ is the associated semigroup, and the final equality holds as a simple consequence of
the fact that $\frac{d}{dt}(P^F_t\varphi_k)=-P^F_t\Delta\varphi_k=-\lambda_k P^F_t\varphi_k$.
Now by \eqref{asum-a}, it holds that $1=P_t^F1$ is in the domain of $\Delta$. A standard argument
thus yields $\Delta 1 = \Delta P_t^F1=-\frac{d}{dt}(P^F_t1)=0$, and so there is no loss of generality in
presupposing that $\lambda_0=0$ and $\varphi_0\equiv1$ in this setting. The only additional assumption
we make on the transition density $(q_t(x,y))_{x,y\in F,t>0}$ is that it is jointly continuous in $(t,x,y)$
(i.e. \eqref{asum-b} holds).

{\lem \label{speclem} Suppose that the operator $\Delta$ has a compact resolvent, so that the above spectral
decomposition holds. If there is a spectral gap, i.e. $\lambda_1>0$, then $(q_t(x,y))_{x,y\in F,t>0}$ converges
to stationarity in an $L^p$ sense (namely \eqref{erg} holds) for any $p\in[1,\infty]$.}
\begin{proof} Recall from (\ref{l2obs}) that $D_2(x,t)^2=q_{2t}(x,x)-1$. Under the assumptions of the lemma,
it follows that
\begin{equation}
D_2(x,t)^2=\sum_{k\geq 1}e^{-2\lambda_kt}\varphi_k(x)^2\rightarrow 0,\label{l2conv}
\end{equation}
as $t\rightarrow\infty$, which completes the proof of the result for $p=2$. To extend this to any $p$, we
first use Cauchy-Schwarz to deduce
\begin{eqnarray*}
(q_t(x,y)-1)^2&=&\left(\sum_{k\geq 1}e^{-\lambda_kt}\varphi_k(x)\varphi_k(y)\right)^2\\
&\leq&
\sum_{k\geq 1}e^{-\lambda_kt}\varphi_k(x)^2\sum_{k\geq 1}e^{-\lambda_kt}\varphi_k(y)^2\\
&=&(q_t(x,x)-1)(q_t(y,y)-1).
\end{eqnarray*}
Consequently, we have that
\begin{eqnarray*}
D_\infty(x,t)^2&=&\sup_{y\in F}(q_t(x,y)-1)^2\\
&\leq & (q_t(x,x)-1)\sup_{y\in F}(q_t(y,y)-1)\\
&\leq & D_2(x,t/2)^2\sup_{y\in F}D_\infty(y,1)\\
\end{eqnarray*}
for any $t\geq 1$, where the second inequality involves an application of the monotonicity property proved as
part of Lemma \ref{dplem}. Now, by
\eqref{asum-b}, the term $\sup_{y\in F}D_\infty(y,1)$ is a finite constant, and so combining the above bound
with (\ref{l2conv}) implies that $D_\infty(x,t)\leq CD_2(x,t/2)\rightarrow 0$ as $t\rightarrow\infty$. The result
for general $p\in[1,\infty]$ is an immediate consequence of this.
\end{proof}

We now give a lemma that explains how to check conditions \eqref{asum-c} and \eqref{asum-d}.

{\lem \label{speclem-22} Suppose that the operator $\Delta$ has a compact resolvent and there is a spectral
gap, then the conditions \eqref{asum-c} and \eqref{asum-d} are automatically satisfied.}
\begin{proof}
Firstly, assume that $q_t(x,y)=0$ for some $x,y\in F$, $t>0$. If $s\in (0,t)$, then the Chapman-Kolmogorov
equations yield $0=q_t(x,y)=\int_Fq_s(x,z)q_{t-s}(z,y)\pi(dz)$.
Since $\pi$ has full support,
using \eqref{asum-b}, it follows that $q_s(x,z)q_{t-s}(z,y)=0$ for every $z\in F$. In particular, $q_s(x,y)q_{t-s}(y,y)=0$.
Noting that $q_{t-s}(y,y)=D_2^2(y,t/2)+1\geq 1$, we deduce that $q_s(x,y)=0$. Now, define a function
$f:(0,\infty)\rightarrow\mathbb{R}_+$ by setting $f(s):=q_s(x,y)$. Letting $(\lambda_i')_{i\geq 0}$ represent
the distinct eigenvalues of $\Delta$, we can write
\[f(s)=\sum_{i\geq 0}a_ie^{-\lambda_i's},\]
where $a_i:=\sum_{j:\lambda_j=\lambda_i'}\varphi_j(x)\varphi_j(y)$. In fact, since Cauchy-Schwarz implies
$\sum_{i\geq 0}|a_ie^{-\lambda_i's}|\leq (q_s(x,x)q_s(y,y))^{1/2}<\infty$, this series converges absolutely
whenever $s\in(0,\infty)$. Thus $f(z):=\sum_{i\geq 0}a_ie^{-\lambda_i'z}$ defines an analytic function on
the whole half-plane $\Re(z)>0$. By our previous observation regarding $q_s(x,y)$, this analytic function is
equal to 0 on the set $(0,t]$, and therefore it must be 0 everywhere on $\Re(z)>0$. However, this contradicts
the fact that $f(t)=q_t(x,y)\rightarrow 1$ as $t\rightarrow\infty$, which was proved in Lemma \ref{speclem}.
Hence, $q_t(x,y)>0$ for every $x,y\in F$, $t>0$.

Secondly, suppose that $q_t(x,\cdot)\equiv 1$ for some $x\in F$ and $t>0$. Then $1=q_t(x,x)=1+
\sum_{i\geq 1}\varphi_i(x)^2e^{-\lambda_i t}$, and so $\varphi_i(x)=0$ for every $i\geq 1$. This implies
that $q_t(x,x)=1$ for every $t>0$. However, by following the proof of (\ref{0lim}), one can deduce that
\[\lim_{t\rightarrow0} (q_t(x,x)-1)=\lim_{t\rightarrow0}D_2^2(x,t/2)\geq\lim_{t\rightarrow0} D_1^2(x,t/2)\geq2,\]
and so the previous conclusion can not hold. Consequently, we have shown that $q_t(x,\cdot)\not\equiv 1$ for any
$x\in F$, $t>0$, as desired.
\end{proof}

To summarize, the above results demonstrate that to verify all the conditions on the transition density that
are required to apply our mixing time convergence results, it will suffice to check that the conservative Hunt
process $X^F$ has a jointly continuous transition density and the corresponding non-negative self-adjoint
operator, $\Delta$, has a compact resolvent and exhibits a spectral gap. As the following corollary explains,
this is a particularly useful observation in the case that the Dirichlet form $(\mathcal{E},\mathcal{F})$
associated with $X^F$ is a resistance form. A precise definition of such an object appears in \cite[Definition 3.1]{KigRes}, for example, but the key property is the finiteness of the corresponding resistance, i.e.
\[R(x,y):=\sup\left\{\frac{|f(x)-f(y)|^2}{\mathcal{E}(f,f)}\: :\:f\in\mathcal{F},\: \mathcal{E}(f,f)>0\right\}\]
is finite for any $x,y\in F$.

{\cor\label{resform} Suppose that $X^F$ is a $\pi$-symmetric Hunt process on $F$ such that the associated
Dirichlet form $(\mathcal{E},\mathcal{F})$ is a resistance form, then \eqref{asum-a}-\eqref{erg} are
automatically satisfied.}
\begin{proof} The fact that $X^F$ is conservative is clear since for a resistance form $1\in  \mathcal{F}$
and $\mathcal{E}(1,1)=0$. That \eqref{asum-b} holds is proved in \cite[Lemma 10.7]{KigRes}. Moreover,
we can check that the non-negative operator corresponding to $(\mathcal{E},\mathcal{F})$ has a compact
resolvent (see \cite[Lemma 9.7]{KigRes} and \cite[Theorem B.1.13]{Kigami}) and exhibits a spectral gap
(this is an easy consequence of the fact that, for a resistance form, $\mathcal{E}(f,f)=0$ if and only if $f$
is constant). Thus, by Lemma \ref{speclem} and Lemma \ref{speclem-22}, the transition density of $X^F$
also satisfies \eqref{asum-c}-\eqref{erg}.
\end{proof}

\section{Examples}\label{examsec}

The mixing time results of the previous sections have many applications. To begin with a particularly simple one, consider $G^N$ to be a discrete $d$-dimensional box of side-length $N$, i.e. vertex set $\{1,2,\dots,N\}^d$ and nearest neighbor connections. By applying classical results about the convergence of the simple random walk on this graph to Brownian motion on $[0,1]^d$ reflected at the boundary, Theorem \ref{main} readily implies that the $L^p$-mixing time of the simple random walk on $\{1,2,\dots,N\}^d$, when rescaled by $N^{-2}$, converges to the $L^p$-mixing time of the limit process for any $p\in [1,\infty]$. A similar result could be proved for the random walk on the discrete torus $(\mathbb{Z}/N\mathbb{Z})^d$. More interestingly, however, as we will now demonstrate, it is possible to apply our main results in a number of examples where the graphs, and sometimes limiting spaces, are random: self-similar fractal graphs with random weights, critical Galton-Watson trees, the critical Erd\H{o}s-R\'{e}nyi random graph, and the range of the random walk in high dimensions. For the second and third of these, we will in the next section go on to describe how the convergence in distribution of mixing times we establish can be applied to relate tail asymptotics for mixing time distributions of the discrete and continuous models.

\subsection{Self-similar fractal graphs with random weights}\label{sssec}

Although the results we have proved apply more generally to self-similar fractal graphs (see below for some further comments on this point), to keep the presentation concise we restrict our attention here to graphs based on the classical Sierpinski gasket, the definition of which we now recall. Suppose $p_1,p_2,p_3$ are the vertices of an equilateral triangle in $\mathbb{R}^2$. Define the similitudes
\[\psi_i(x):=p_i+\frac{z-p_i}{2},\hspace{20pt}i=1,2,3.\]
Since $(\psi_i)_{i=1}^3$ is a family of contraction maps, there exists a unique non-empty compact set $F$ such that $F=\cup_{i=1}^3 \psi_i(F)$ -- this is the Sierpinski gasket. We will suppose $d_F$ is the intrinsic shortest path metric on $F$ defined in \cite{Kigamimetric}, and note that this induces the same topology as the Euclidean metric. Moreover, we suppose $\pi$ is the $(\ln 3)/(\ln 2)$-Hausdorff measure on $F$ with respect to the Euclidean metric, normalized to be a probability measure. This measure is non-atomic, has full support and satisfies $\pi(\partial B(x,r))=0$ for every $x\in F$, $r>0$ (see \cite[Lemma 25]{CHLLT}).

We now define a sequence of graphs $(G^N)_{N\geq 0}$ by setting
\[V(G^N):=\bigcup_{i_1,\dots ,i_N=1}^{3}\psi_{i_1\dots i_N}(V_0),\]
where $V_0:=\{p_1,p_2,p_3\}$ and $\psi_{i_1\dots i_n}:=\psi_{i_1}\circ\dots\circ\psi_{i_n}$, and
\[E(G^N):=\left\{\{\psi_{i_1\dots i_N}(x),\psi_{i_1\dots i_N}(y)\}:\:x,y\in V_0,\:x\neq y,\:i_1,\dots,i_N\in\{1,2,3\}\right\}.\]
We set $d_{G^N}:=d_F|_{V(G^N)\times V(G^N)}$, so that $(V(G^N),d_{G^N})$ converges to $(F,d_F)$ with respect to the Hausdorff distance between compact subsets of $F$. Weights $(\mu^N_e)_{e\in E(G^N),N\geq 0}$ will be selected independently at random from a common distribution, which we assume is supported on an interval $[c_1,c_2]$, where $0<c_1\leq c_2<\infty$. By the procedure described in the introduction, we define from these weights a sequence of random measures $(\pi^N)_{N\geq 0}$ on the vertex sets of our graphs in the sequence $(G^N)_{N\geq 0}$. That $\pi^N$ weakly converges to $\pi$ as Borel probability measures on $F$, almost-surely, can be checked by applying \cite[Lemma 26]{CHLLT}.

To describe the scaling limit of the random walks associated with the random weights $\mu^N$, we appeal to the homogenization result of \cite{homog}. To describe this, we first introduce the Dirichlet form associated with the walk on the level $N$ graph by setting, for $f\in \mathbb{R}^{V(G^N)}$,
\begin{equation}\label{levelN}
\mathcal{E}^N(f,f):=\sum_{i_1,\dots,i_N=1}^3\sum_{x,y\in V_0, x\neq y}\mu^N_{\psi_{i_1\dots i_N}(x)\psi_{i_1\dots i_N}(y)}\left(f(\psi_{i_1\dots i_N}(x))-f(\psi_{i_1\dots i_N}(y))\right)^2.
\end{equation}
Let $\Lambda^N=(\Lambda^N_{xy})_{x,y\in V_0,x\neq y}$ be the collection of weights such that the associated random walk on $G^0$ is the trace of $X^{G^N}$ onto $V_0$. It then follows from \cite[Theorem 3.4]{homog} that there exists a deterministic constant $C\in(0,\infty)$ such that
\[\lim_{n\rightarrow\infty}\mathbf{E}\left|\left(\frac{5}{3}\right)^N\Lambda^N_{xy}-C\right|=0,\]
for any $x,y\in V_0$, $x\neq y$. Now, suppose $\mathcal{E}^N_{C}$ is a quadratic form on $\mathbb{R}^{V(G^N)}$ which satisfies (\ref{levelN}) with $\mu^N_{\psi_{i_1\dots i_N}(x)\psi_{i_1\dots i_N}(y)}$ replaced by $C$ in each summand, then define
\[\mathcal{E}(f,f)=\lim_{N\rightarrow\infty}\left(\frac{5}{3}\right)^N\mathcal{E}^N_{C}(f|_{V(G^N)},f|_{V(G^N)})\]
for $f\in \mathcal{F}$, where $\mathcal{F}$ is the subset of $C(F,\mathbb{R})$ such that the right-hand side above exists and is finite. It is known that $(\mathcal{E},\mathcal{F})$ is a local, regular Dirichlet form on $L^2(F,\pi)$, which is also a resistance form (see \cite{Kigami}, for example). Thus, by Corollary \ref{resform}, the associated $\pi$-symmetric diffusion $X^F$, which (modulo the scaling constant $C$) is known as Brownian motion on the Sierpinski gasket, satisfies \eqref{asum-a}-\eqref{erg}.

For the case of unbounded fractal graphs, a probabilistic version of \eqref{CLT-dist} was proved as \cite[Proposition 30(i)]{CHLLT} by applying the homogenization result for processes of \cite{KumKus} (cf. \cite{homog}). Since the Sierpinski gasket is a finitely ramified fractal, it is a relatively straightforward technical exercise to adapt this result to the compact case by considering a decomposition of the sample paths of the relevant processes into segments started at one of the outer corners of the gasket and stopped upon hitting another.

To expand on this, we will explain how to prove a version of \cite[Theorem 3.6]{KumKus} in our setting. (Note that our $X^{G^N}$ is a discrete time Markov chain with $\pi^N$ as the invariant measure, whereas in \cite{KumKus} it was the continuous-time Markov chains with normalized counting measure as the invariant measure that were studied. However, since both measures are comparable and they converge to $\pi$ almost-surely, this difference can be easily resolved.) Recall $p_1$ and $p_2$ are two distinct elements of $V_0$. Let $\sigma_{p_1}^{(0)}(X^{G^N})$ be the first hitting time of $p_1$ by $X^{G^N}$, and for each
$i\in \mathbb{N}$, define inductively
\begin{eqnarray*}
\sigma_{p_2}^{(i)}(X^{G^N})&:=&\inf\left\{m\ge \sigma_{p_1}^{(i-1)}(X^{G^N}): X^{G^N}_m=p_2\right\},\\
\sigma_{p_1}^{(i)}(X^{G^N})&:=&\inf\left\{m\ge \sigma_{p_2}^{(i)}(X^{G^N}): X^{G^N}_m=p_1\right\}.
\end{eqnarray*}
Then, we can write, for continuous $f:F\rightarrow \mathbb{R}$,
\begin{eqnarray}
\lefteqn{\mathbf{E}^{G^N}_{x_N}\left[f\left(X^{G^N}_{5^Nt}\right)\right]}\nonumber\\
&=&{\mathbf E}^{G^N}_{x_N}\left[f\left(X^{G^N}_{5^Nt}\right): {5^Nt}<\sigma_{p_1}^{(0)}\right]+
{\mathbf E}^{G^N}_{x_N}\left[f\left(X^{G^N}_{5^Nt}\right):\sigma_{p_1}^{(0)}\le {5^Nt}<\sigma_{p_2}^{(1)}\right]~~~~\label{eq:gsjwbw0}\\
&&+\sum_{i=1}^\infty {\mathbf E}^{G^N}_{x_N}\left[f\left(X^{G^N}_{5^Nt}\right):\sigma_{p_2}^{(i)}\le {5^Nt}<\sigma_{p_1}^{(i)}\right]~~~~\label{eq:gsjwbw}\\
&&+\sum_{i=1}^\infty
{\mathbf E}^{G^N}_{x_N}\left[f\left(X^{G^N}_{5^Nt}\right):\sigma_{p_1}^{(i)}\le {5^Nt}<\sigma_{p_2}^{(i+1)}\right],~~~~\label{eq:gsjwbw2}
\end{eqnarray}
where $x_N\in V(G^N)$ converges to $x\in F$, say. The first summand in the right hand side of \eqref{eq:gsjwbw0} can be written in terms of the process $X^{G^N}$ killed at $p_1$, and so by tracing the proof of \cite[Proposition 30(i)]{CHLLT} line by line, we can check it converges to the corresponding expectation involving $X^F$ killed on hitting $p_1$. Similarly, the second summand in \eqref{eq:gsjwbw0} can be written as
\begin{eqnarray*}
\lefteqn{{\mathbf E}^{G^N}_{x_N}[f(X^{G^N}_{5^Nt}):\sigma_{p_1}^{(0)}\le {5^Nt}<\sigma_{p_2}^{(1)}]}\\
&=&
{\mathbf E}^{G^N}_{x_N}[1_{\{\sigma_{p_1}^{(0)}\le {5^Nt}\}}{\mathbf E}^{G^N}_{p_1}
[f(X^{G^N}_{{5^Nt}-\sigma_{p_1}^{(0)}})1_{\{{5^Nt}-\sigma_{p_1}^{(0)}<\sigma_{p_2}^{(1)}\circ \,\theta_{\sigma_{p_1}^{(0)}}\}}|{\cal F}_{\sigma_{p_1}^{(0)}}]],
\end{eqnarray*}
where $\theta$ is the shift map. Given $\sigma_{p_1}^{(0)}=s$, the strong Markov property allows us to write
${\mathbf E}^{G^N}_{p_1}
[f(X^{G^N}_{{5^Nt}-s})1_{\{{5^Nt}-s<\sigma_{p_2}^{(1)}\}}]$ in terms of the process started at $p_1$ and killed on hitting $p_2$, independently of
the distribution of $\sigma_{p_1}^{(0)}$. Thus the second term in the right hand side of \eqref{eq:gsjwbw0} converges to
${\mathbf E}^F_{x}[f(X^{F}_t):\sigma_{p_1}^{(0)}(X^F)\le t<\sigma_{p_2}^{(1)}(X^F)]$. We can prove convergence of the rest of the terms similarly. Moreover, by applying the estimate for the exit time of the random walks from balls stated as part of \cite[Lemma 27]{CHLLT}, for example, it is straightforward to check that there exists a $t_0>0$ such that $\mathbf{P}^{G^N}_{p_1}(\sigma_{p_2}^{(1)}\leq 5^N t_0)$ and $\mathbf{P}^{G^N}_{p_2}(\sigma_{p_1}^{(0)}\leq 5^N t_0)$ are both bounded above by $1/2$, uniformly in $N$. As a consequence of this, one can show that the terms in the sums at \eqref{eq:gsjwbw} and \eqref{eq:gsjwbw2} decay exponentially, uniformly in $N$, and hence that the right hand side of \eqref{eq:gsjwbw} converges to ${\mathbf E}^F_{x}[f(X^{F}_t)]$ as $N\to\infty$.
Convergence of the finite dimensional distributions can be shown similarly and we obtain the desired version of \cite[Theorem 3.6]{KumKus}.

Finally, a probabilistic version of the tightness condition of (\ref{tightcond}) is easily checked by applying (a probabilistic version of) Lemma \ref{tightlem}, using known resistance estimates for nested fractals (cf. \cite[Proposition 30(ii)]{CHLLT}), and so Assumption 1 holds in probability due to Proposition \ref{local}. Thus we are able to apply Theorem \ref{main} to deduce the following.

{\thm If $t_{\rm mix}(G^N)$ is the mixing time of the random walk on the level $N$ approximation to the Sierpinski gasket equipped with uniformly bounded, independently and identically distributed random weights, then
\[5^{-N}t_{\rm mix}(G^N)\rightarrow t_{\rm mix}(F)\]
in probability, where $t_{\rm mix}(F)$ is the mixing time of the diffusion $X^F$.}
\bigskip

Let us remark that the same argument will yield at least two generalizations of this theorem. Firstly, it is not necessary for the weights to be independent and identically distributed, but rather it will be sufficient for them only to be `cell independent', i.e. each collection $(\mu_{\psi_{i_1\dots i_N}(x)\psi_{i_1\dots i_N}(y)}^N)_{x,y\in V_0,x\neq y}$ is independent and identically distributed as $(\mu_{xy})_{x,y\in V_0,x\neq y}$. (We note that without a symmetry condition, though, the limiting diffusion will no longer be guaranteed to be the Brownian motion on the Sierpinski gasket.) Secondly, the Sierpinski gasket is just one example of a nested fractal. Identical arguments could be applied to obtain corresponding mixing time results for sequences of graphs based on any of the highly-symmetric fractals that come from this class (since the key references \cite{CHLLT}, \cite{homog} and \cite{KumKus} all incorporate nested fractals already).

Finally, variations on the above mixing time convergence result can also be established for examples along the lines of those appearing in \cite[Sections 7.4 and 7.5]{CHLLT}. These include: an almost-sure statement for Vicsek set-type graphs (which complements the mixing time bounds for deterministic versions of these graphs proved in \cite{GMT}); a convergence of mixing times for deterministic Sierpinski carpet graphs; and a subsequential limit for Sierpinski carpets with random weights. Since many of the ideas needed for these applications are similar to those discussed above, we omit the details.

\subsection{Critical Galton-Watson trees}\label{treesec}

The connection between critical Galton-Watson processes and $\alpha$-stable trees is now well-known, and so we will be brief in introducing it. Let $\xi$ be a mean 1 random variable whose distribution is aperiodic (not supported on a sub-lattice of $\mathbb{Z}$). Furthermore, suppose that $\xi$ is in the domain of attraction of a stable law with index $\alpha\in(1,2)$, by which we mean that there exists a sequence $a_N\rightarrow\infty$ such that
\begin{equation}\label{stable}
\frac{\xi[N]-N}{a_N}\rightarrow\Xi,
\end{equation}
in distribution, where $\xi[N]$ is the sum of $N$ independent copies of $\xi$ and the limit random variable satisfies $\mathbf{E}(e^{-\lambda\Xi})=e^{-\lambda^\alpha}$ for $\lambda>0$. If $\mathcal{T}_N$ is a Galton-Watson tree with offspring distribution $\xi$ conditioned to have total progeny $N$, then it is the case that
\begin{equation}\label{treeconv}
N^{-1}a_N\mathcal{T}_N\rightarrow \mathcal{T}^{(\alpha)},
\end{equation}
in distribution with respect to the Gromov-Hausdorff distance between compact metric spaces, where $\mathcal{T}^{(\alpha)}$ is an $\alpha$-stable tree normalized to have total mass equal to 1 (see \cite[Theorem 4.3]{rrt}, which is a corollary of a result originally proved in \cite{Duqap}). Note that the left-hand side here is shorthand for the metric space $(V(\mathcal{T}_N),N^{-1}a_Nd_{\mathcal{T}_N})$, where $V(\mathcal{T}_N)$ is the vertex set of $\mathcal{T}_N$ and $d_{\mathcal{T}_N}$ is the shortest path graph distance on this set.

The $\alpha$-stable tree $\mathcal{T}^{(\alpha)}$ is almost-surely a compact metric space. Moreover, there is a natural non-atomic probability measure upon it, $\pi^{(\alpha)}$ say, which has full support, and appears as the limit of the uniform measure on the approximating graph trees. Usefully, we can decompose this measure in terms of a collection of measures of level sets of the tree. More specifically, in the construction of the $\alpha$-stable tree from an excursion we can naturally choose a root $\rho\in \mathcal{T}^{(\alpha)}$. We define $\mathcal{T}^{(\alpha)}(r):=\{x\in \mathcal{T}^{(\alpha)}:d_{\mathcal{T}^{(\alpha)}}(\rho,x)=r\}$ to be the collection of vertices at height $r$ above this vertex. For almost-every realization of $\mathcal{T}^{(\alpha)}$, there then exists a c\`{a}dl\`{a}g sequence of finite measures on $\mathcal{T}^{(\alpha)}$, $(\ell^r)_{r>0}$, such that $\ell^r$ is supported on $\mathcal{T}^{(\alpha)}(r)$ for each $r$ and
\[\pi^{(\alpha)}=\int_0^\infty \ell^r dr\]
(see \cite[Section 4.2]{LegallDuquesne}). Clearly this implies that $\pi^{(\alpha)}(\partial B_{\mathcal{T}^{(\alpha)}}(\rho,r))=0$ for every $r>0$, for almost-every realization of $\mathcal{T}^{(\alpha)}$. Since $\alpha$-stable trees satisfy a root-invariance property (see \cite[Theorem 4.8]{LegallDuquesne}), one can easily extend this result to hold for $\pi^{(\alpha)}$-a.e. $x\in \mathcal{T}^{(\alpha)}$. Although this is not quite the assumption of the introduction that $\pi^{(\alpha)}(\partial B_{\mathcal{T}^{(\alpha)}}(x,r))=0$ for every $x\in \mathcal{T}^{(\alpha)}$, $r>0$, by a minor tweak of the proof of Proposition \ref{point}, we are still able to apply our mixing time convergence results in the same way.

Upon almost-every realization of the metric measure space $(\mathcal{T}^{(\alpha)}, \pi^{(\alpha)})$, it is possible to define a corresponding Brownian motion $X^{(\alpha)}$ (to do this, apply \cite[Theorem 5.4]{Kigamidendrite}, in the way described in \cite[Section 2.2]{Croydoncbp}). This is a conservative $\pi^{(\alpha)}$-symmetric Hunt process, and the associated Dirichlet form $(\mathcal{E}^{(\alpha)},\mathcal{F}^{(\alpha)})$ is actually a resistance form. Thus we can again apply Corollary \ref{resform} to confirm that
\eqref{asum-a}-\eqref{erg} hold
for some corresponding transition density, $q^{(\alpha)}$ say. Now, in \cite{Croydoninf}, it was demonstrated that if $\mathbf{P}^{\mathcal{T}_N}_{\rho^{N}}$ is the law of the random walk on $\mathcal{T}_N$ started from its root (original ancestor) $\rho_N$ and $\pi^N$ is its stationary probability measure, then, after embedding all the objects into an underlying Banach space in a suitably nice way, the conclusion of (\ref{treeconv}) can be extended to the distributional convergence of
\begin{equation}\nonumber
\left(N^{-1}a_N\mathcal{T}_N,\pi^N(Na_N^{-1}\cdot), \mathbf{P}^{\mathcal{T}_N}_{\rho^{N}}\left(\left(N^{-1}a_NX^{\mathcal{T}_N}_{\lfloor N^2a_N^{-1} t\rfloor}\right)_{t\in [0,1]}\in \cdot\right)\right)
\end{equation}
to $(\mathcal{T}^{(\alpha)},\pi^{(\alpha)},\mathbf{P}_\rho^{(\alpha)})$, where $\mathbf{P}_\rho^{(\alpha)}$ is the law of $X^{(\alpha)}$ started from $\rho$. By applying the fixed starting point version of the local limit result of Proposition \ref{local} (cf. \cite[Theorem 1]{CHLLT}), similarly to the argument of \cite[Section 7.2]{CHLLT}, for the Brownian continuum random tree, which corresponds to the case $\alpha=2$, one can obtain from this a distributional version of Assumption \ref{assu3}. (The tightness condition of (\ref{tightcond}) is easily checked by applying Lemma \ref{tightlem}.)

{\lem \label{treedist}For any compact interval $I\subset (0,\infty)$,
\[{\left(\left(V(\mathcal{T}_N),N^{-1}a_Nd_{\mathcal{T}^N},\rho^N\right), \pi^N, \left(q^N_{N^{2}a_N^{-1}t}(\rho^N,x)\right)_{x\in V(\mathcal{T}_N),t\in I}\right)}\]
converges in distribution to $((\mathcal{T}^{(\alpha)},d_{\mathcal{T}^{(\alpha)}},\rho),\pi^{(\alpha)},(q^{(\alpha)}_t(\rho,x))_{x\in \mathcal{T}^{(\alpha)},t\in I})$ in a spectral pointed Gromov-Hausdorff sense.}
\bigskip

Consequently, since the space in which the above convergence in distribution occurs is separable, we can use a Skorohod coupling argument to deduce from this and Theorem \ref{pointthm} the following mixing time convergence result. We remark that the $\sqrt 2$ that appears in the finite variance result is simply an artefact of the particular scaling we have described here, and could alternatively have been absorbed in the scaling of metrics.

{\thm\label{treethm} Fix $p\in [1,\infty]$. If $t_{\rm mix}^{p}(\rho^N)$ is the $L^p$-mixing time of the random walk on $\mathcal{T}_N$ started from its root $\rho^N$, then
\[N^{-2}a_Nt_{\rm mix}^{p}(\rho^N) {\rightarrow} t_{\rm mix}^{p}(\rho),\]
in distribution, where $t_{\rm mix}^{p}(\rho)\in(0,\infty)$ is the  $L^p$-mixing time of the Brownian motion on $\mathcal{T}^{(\alpha)}$ started from $\rho$. In particular, in the case when the offspring distribution has finite variance $\sigma$, it is the case that
\[\frac{\sigma}{\sqrt{2}} N^{-3/2}t_{\rm mix}^{p}(\rho^N){\rightarrow} t_{\rm mix}^{p}(\rho),\]
in distribution.}
\bigskip

\begin{rem}\label{rem:ext}
{\rm We note that it was only for convenience that the convergence of the random walks on the trees
$\mathcal{T}_N$, $N\geq 1$, to the Brownian motion on $\mathcal{T}^{(\alpha)}$ was proved
from a single starting point in \cite{Croydoninf}. We do not anticipate any significant problems in
extending this result to hold for arbitrary starting points. Indeed, the first step would be to make
the obvious adaptations to the proof of \cite[Lemma 4.2]{Croydoninf} to extend the result, which
demonstrates convergence of simple random walks (and related additive functionals) on subtrees
of $\mathcal{T}_N$ consisting of a finite number of branch segments to the corresponding continuous
objects, from the case when all the random walks start from the root to an arbitrary starting point
version. An argument identical to the remainder of \cite[Section 4]{Croydoninf} could then be used
to obtain the convergence of simple random walks on the whole trees, at least in the case when the
starting point of the diffusion is in one of the finite subtrees considered. Since the union of
the finite subtrees is dense in the limiting space, we could subsequently use the heat kernel
continuity properties to obtain the non-pointed spectral Gromov-Hausdorff version of Lemma~\ref{treedist}.
However, we do not pursue this approach here as it would require a substantial amount of space and new
notation that is not relevant to the main ideas of this article. Were it to be checked, Theorem~\ref{main}
would imply, for any $p\in[1,\infty]$, the distributional convergence of $t_{\rm mix}^{p}(\mathcal{T}_N)$,
the $L^p$-mixing time of the random walk on $\mathcal{T}_N$, when rescaled appropriately,
to $t_{\rm mix}^{p}(\mathcal{T}^{(\alpha)})\in(0,\infty)$, the  $L^p$-mixing time of the
Brownian motion on $\mathcal{T}^{(\alpha)}$.}
\end{rem}

\subsection{Critical Erd\H{o}s-R\'{e}nyi random graph}\label{ersec}

Closely related to the random trees of the previous section is the Erd\H{o}s-R\'{e}nyi random graph at criticality. In particular, let $G(N,p)$ be the random graph in which every edge of the complete graph on $N$ labeled vertices $\{1,\dots, N\}$ is present with probability $p$ independently of the other edges. Supposing $p=N^{-1}+\lambda N^{-4/3}$ for some $\lambda\in \mathbb{R}$, so that we are in the so-called critical window, it is known that the largest connected component $\mathcal{C}^N$,
equipped with its shortest path graph metric $d_{\mathcal{C}^N}$, satisfies
\[\left(V(\mathcal{C}^N),N^{-1/3}d_{\mathcal{C}^N}\right)\rightarrow \left(\mathcal{M},d_\mathcal{M}\right)\]
in distribution, again with respect to the Gromov-Hausdorff distance between compact metric spaces, where $(\mathcal{M},d_\mathcal{M})$ is a random compact metric space \cite{ABG}.
(In fact, this and all the results given in this subsection hold for a family of $i$-th largest connected components for all $i\in \mathbb{N}$.
For notational simplicity, we only discuss the largest connected component $\mathcal{C}^N$.)
Moreover, in \cite{Croydoncrg}, it was shown that the associated random walks started from a root vertex $\rho^N$ satisfy a distributional convergence result of the form
\[\left(N^{-1/3}X^{\mathcal{C}^N}_{{\lfloor Nt \rfloor}}\right)_{t\geq 0}\rightarrow \left(X^{\mathcal{M}}_t\right)_{t\geq 0},\]
where $X^\mathcal{M}$ is a diffusion on the space $\mathcal{M}$ started from a distinguished vertex $\rho\in\mathcal{M}$. Although the invariant probability measures of the random walks, $\pi^N$ say, were not considered in \cite{Croydoncrg}, it is not difficult to extend this result to include them since the hard work regarding their convergence has already been completed (see \cite[Lemma 6.3]{Croydoncrg}, in particular). Hence, by again applying the fixed starting point version of the local limit result of Proposition \ref{local} (using Lemma \ref{tightlem} again to deduce the relevant tightness condition), we are able to obtain the analogue of Lemma \ref{treedist} in this setting.

{\lem For any compact interval $I\subset (0,\infty)$,
\[{\left(\left(V(\mathcal{C}^N),N^{-1/3}d_{\mathcal{C}^N},\rho^N\right), \pi^N, \left(q^N_{Nt}(\rho^N,x)\right)_{x\in V(\mathcal{T}_N),t\in I}\right)},\]
converges in distribution to $((\mathcal{M},d_{\mathcal{M}},\rho),\pi^\mathcal{M},(q^\mathcal{M}_t(\rho,x))_{x\in \mathcal{M},t\in I})$, where $\pi^\mathcal{M}$ is the invariant probability measure of $X^\mathcal{M}$ and $(q^\mathcal{M}_t(x,y))_{x,y\in \mathcal{M},t>0}$ is its transition density with respect to this measure, in a spectral pointed Gromov-Hausdorff sense.}
\bigskip

In order to proceed as above, we must of course check that $\pi^\mathcal{M}$ and $q^\mathcal{M}$ satisfy a number of technical conditions. To do this, first observe that a typical realization of $\mathcal{M}$ looks like a (rescaled) typical realization of the Brownian continuum random tree $\mathcal{T}^{(2)}$ glued together at a finite number of pairs of points \cite{ABG}. Since $\pi^\mathcal{M}$ can be considered as the image of the canonical measure $\pi^{(2)}$ on $\mathcal{T}^{(2)}$ under this gluing map, it is elementary to obtain from the statements of the previous section regarding $\pi^{(2)}$ that $\pi^\mathcal{M}$ is almost-surely non-atomic, has full support and satisfies $\pi^\mathcal{M}(\partial B_\mathcal{M}(x,r))=0$ for $\pi^\mathcal{M}$-a.e. $x\in \mathcal{M}$ and every $r>0$, as desired. For $q^\mathcal{M}$, we simply observe that because the Dirichlet form corresponding to $X^\mathcal{M}$ is a resistance form (\cite[Proposition 2.1]{Croydoncrg}), we can once again apply Corollary \ref{resform} to establish conditions \eqref{asum-a}-\eqref{erg}.

Given these results, pointwise mixing time convergence follows from Theorem \ref{pointthm}.

{\thm \label{c1nmix} Fix $p\in [1,\infty]$. If $t_{\rm mix}^{p}(\rho^N)$ is the $L^p$-mixing time of the random walk on $\mathcal{C}^N$ started from its root $\rho^N$, then
\[N^{-1}t_{\rm mix}^{p}(\rho^N) {\rightarrow} t_{\rm mix}^{p}(\rho),\]
in distribution, where $t_{\rm mix}^{p}(\rho)\in(0,\infty)$ is the  $L^p$-mixing time of the Brownian motion on $\mathcal{M}$ started from $\rho$.}
\bigskip

\begin{rem}{\rm
As discussed in Remark~\ref{rem:ext}, we do not expect any major barriers in extending the above result to
arbitrary starting points. The first task in doing this would be to adapt the convergence result proved in
\cite{Croydoncrg} regarding the convergence of simple random walks on subgraphs of $\mathcal{C}_1^n$
formed of a finite number of line segments (\cite[Lemma 6.4]{Croydoncrg}) to arbitrary starting points.
One could then extend this to obtain the desired convergence result for simple random walks on the entire
space using ideas from \cite[Section 7]{Croydoncrg} and heat kernel continuity. It would also be necessary
to introduce a new Gromov-Hausdorff-type topology to state the result, as the one used in \cite{Croydoncrg}
is only suitable for the pointed case. Again, we suspect taking these steps will simply be a lengthy
technical exercise, and choose not to follow them through here. We do though reasonably expect that
$t_{\rm mix}^{p}(\mathcal{C}^N)$, the $L^p$-mixing time of the random walk on $\mathcal{C}^N$,
when rescaled appropriately, converges in distribution to $t_{\rm mix}^{p}(\mathcal{M})\in(0,\infty)$,
the  $L^p$-mixing time of the Brownian motion on $\mathcal{M}$, for any $p\in[1,\infty]$.}
\end{rem}

\subsection{Random walk on range of random walk in high dimensions}\label{rrwsec}

Let $S=(S_n)_{n\geq 0}$ be the simple random walk on $\mathbb{Z}^d$ started from 0, built on an underlying probability space with probability measure $\mathbf{P}$, and define the range of $S$ up to time $N$ to be the graph ${G}^N$ with vertex set
\begin{equation}\label{rangev}
V({G}^N):=\left\{S_n:0\leq n\leq N\right\},
\end{equation}
and edge set
\begin{equation}\label{rangee}
E({G}^N):=\left\{\{S_{n-1},S_{n}\}:1\leq n\leq N\right\}.
\end{equation}
In this section, we will explain how to prove that if $d\geq 5$, which is an assumption henceforth, then the mixing times of the sequence of graphs $(G^N)_{N \geq 1}$ grows asymptotically as $cN^2$, $\mathbf{P}$-a.s., where $c$ is a deterministic constant. Since doing this primarily depends on making relatively simple adaptations of the high-dimensional scaling limit result of \cite{rwrrw} for the random walk on the entire range of $S$ (i.e. the $N=\infty$ case) to the finite length setting, we will be brief with the details.

First, suppose that ${S}=({S}_n)_{n\in \mathbb{Z}}$ is a two-sided extension of $({S}_n)_{n\geq0}$ such that $({S}_{-n})_{n\geq 0}$ is an independent copy of $({S}_n)_{n\geq0}$. The set of cut-times for this process, \[{\mathcal{T}}:=\left\{n:{S}_{(-\infty,n]}\cap{S}_{[n+1,\infty)}=\emptyset\right\},\]
is known to be infinite $\mathbf{P}$-a.s. (\cite{ET}). Thus we can write ${\mathcal{T}}=\{T_n:n\in \mathbb{Z}\}$, where $\dots T_{-1}<T_0\leq0<T_1<T_2<\dots$. The corresponding set of cut-points is given by ${C}:=\{C_n:n\in\mathbb{Z}\}$, where $C_n:=S_{T_n}$. For these objects, an ergodicity argument can be applied to obtain that, $\mathbf{P}$-a.s., as $|n|\rightarrow\infty$,
\begin{equation}\label{taud}
\frac{T_n}{n}\rightarrow \tau(d):={\mathbf{E}}(T_1| 0\in{\mathcal{T}})\in[1,\infty),
\end{equation}
\[\frac{d_{{{G}}}(0,C_n)}{|n|}\rightarrow \delta(d):={\mathbf{E}}(d_{{{G}}}(0,C_1)| 0\in\mathcal{T})\in[1,\infty),\]
where $d_{{{G}}}$ is the shortest path graph distance on the range ${G}$ of the entire two-sided walk ${S}$, which is defined analogously to (\ref{rangev}) and (\ref{rangee}). In particular, see \cite[Lemma 2.2]{rwrrw}, for a proof of the same convergence statements under the measure $\mathbf{P}(\cdot|0\in{\mathcal{T}})$, and note that the conditioning can be removed by using the relationship between  $\mathbf{P}$ and $\mathbf{P}(\cdot|0\in{\mathcal{T}})$ described in \cite[Lemma 2.1]{rwrrw}. Given these results, it is an elementary exercise to check that the metric space
$(V(G^N),\tau(d)\delta(d)^{-1}N^{-1}d_{G^N})$, where $d_{G^N}$ is the shortest path graph distance on $G^N$, converges $\mathbf{P}$-a.s. with respect to the Gromov-Hausdorff distance to the interval $[0,1]$ equipped with the Euclidean metric. Moreover, the same ideas readily yield an extension of this result to a spectral Gromov-Hausdorff one including that $\pi^N$, the invariant measure of the associated simple random walk, converges to Lebesgue measure on $[0,1]$.

Now, for a fixed realization of ${G}$, let ${X}=({X}_n)_{n\geq0}$ be the simple random walk on ${G}$ started from 0. Define the hitting times by ${X}$ of the set of cut-points ${\mathcal{C}}$ by $H_0:=\min\{m\geq 0:{X}_m\in{\mathcal{C}}\}$, and, for $n\geq 1$, $H_n:=\min\{m> H_{n-1}:{X}_m\in{\mathcal{C}}\}$. We use these times to define a useful indexing process $Z=(Z_n)_{n\geq0}$ taking values in $\mathbb{Z}$. In particular, if $n<H_0$, define $Z_n$ to be the unique $k\in\mathbb{Z}$ such that ${X}_{H_0}=C_k$. Similarly, if $n\in[H_{m-1},H_{m})$ for some $m\geq 1$, then define $Z_n$ to be the unique $k\in\mathbb{Z}$ such that ${X}_{H_{m}}=C_k$. Noting that this definition precisely coincides with the definition of $Z$ in \cite{rwrrw}, from Lemma 3.5 of that article we have that: for $\mathbf{P}$-a.e. realization of ${G}$,
\begin{equation}\label{zconv}
\left(N^{-1}\tau(d)Z_{\lfloor tN^2\rfloor}\right)_{t\geq0}\rightarrow (B_{t\kappa_2(d)})_{t\geq0},
\end{equation}
in distribution, where $(B_t)_{t\geq0}$ is a standard Brownian motion on $\mathbb{R}$ started from 0, and $\kappa_2(d)\in(0,\infty)$ is the deterministic constant defined in \cite{rwrrw}. To deduce from (\ref{zconv}) the following scaling limit for $X^N$, the simple random walk on $G^N$, we proceed via a time-change argument that is essentially a reworking of parts of \cite[Section 3]{rwrrw}.

{\lem For $\mathbf{P}$-a.e. realization of ${S}$, if $X^N$ is started from 0, then
\[\left(\tau(d)\delta(d)^{-1}N^{-1}d_{G^N}\left(0,X^N_{\lfloor \kappa_2(d)^{-1}N^2t\rfloor}\right)\right)_{t\geq 0}\rightarrow \left(B^{[0,1]}_t\right)_{t\geq 0},\]
in distribution, where $B^{[0,1]}=(B^{[0,1]}_t)_{t\geq 0}$ is Brownian motion on $[0,1]$ started at 0 and reflected at the boundary.}
\begin{proof} The following proof can be applied to any typical realization of ${S}$. To begin with, define a process $(A^{Z,N}_n)_{n\geq0}$ by setting
\[A^{Z,N}_n:=\sum_{m=0}^{n-1}\mathbf{1}_{\{Z_m\in [0,T_N^{-1}]\}},\]
where $T_N^{-1}:=\max\{n:T_n\leq N\}$. From (\ref{taud}), we have that $T_N^{-1}\sim\tau(d)^{-1}N$. Combining this observation with (\ref{zconv}), one can check that, simultaneously with (\ref{zconv}), $(N^{-2}A^N_{\lfloor tN^2\rfloor})_{t\geq0}$ converges in distribution to $(\kappa_2(d)^{-1}A^B_{t\kappa_2(d)})_{t\geq0}$, where
\[A^B_t:=\int_{0}^t\mathbf{1}_{\{B_{s}\in[0,1]\}}ds\]
(cf. \cite[Lemma 3.5]{rwrrw}).

We now apply the above result to establish a scaling limit for the process ${X}$ observed on the vertex set $V(\tilde{G}^N):=\{S_n:T_1\leq n\leq T_{N}^{-1}\}$. Specifically, set
\[{A}^{N}_n:=\sum_{m=0}^{n-1}\mathbf{1}_{\{{X}_m,{X}_{m+1}\in V(\tilde{G}^N)\}}.\]
Similarly to the proof of \cite[Lemma 3.6]{rwrrw}, one can check that  \[\sup_{0\leq m\leq n}\left|{A}^{N}_m-A^{Z,N}_m\right|\leq\sum_{m=0}^{n}\mathbf{1}_{\{Z_m\in[0,1,2]\cup[T_N^{-1}-2,T_N^{-1}-1,T_N^{-1}]\}}.\]
It is therefore a simple consequence of (\ref{zconv}) that $N^{-2}\sup_{0\leq m\leq TN^2}\left|{A}^{N}_m-A^{Z,N}_m\right|$ converges to 0 in probability as $N\rightarrow\infty$ for any $T\in(0,\infty)$. Since we know from equation (16) of \cite{rwrrw} that
\[N^{-1}\sup_{0\leq m\leq TN^2}\left|d_{{G}}\left(0,{X}_{m}\right)-\delta(d)Z_m\right|\]
also converges to 0 in probability, we readily obtain
\begin{equation}\label{tildeconv}
\left(\tau(d)\delta(d)^{-1}N^{-1}d_{{G}}\left(0,\tilde{X}^N_{\lfloor N^2t\rfloor}\right)\right)_{t\geq 0}\rightarrow \left(B^{[0,1]}_{\kappa_2(d)t}\right)_{t\geq 0},
\end{equation}
in distribution, where $\tilde{X}^N=(\tilde{X}^{N}_n)_{n\geq0}$ is the random walk ${X}$ observed on $V(\tilde{G}^N)$ -- this is defined precisely by setting $\tilde{X}^N_n:={X}_{{\alpha}^N(n)}$, where ${\alpha}^N(n):=\max\{{A}_m^N\leq n\}$. We remark that the particular limit process $B^{[0,1]}$ arises as a consequence of the fact that $(B_{\alpha^B(t)})_{t\geq0}$, where $\alpha^B$ is the right-continuous inverse of $A^B$, has exactly the distribution of $B^{[0,1]}$.

Finally, since the process $\tilde{X}^N$ is identical in law to the simple random walk ${X}^N$ observed on $V(\tilde{G}^N)$, to replace $\tilde{X}^N$ by $X^N$ in (\ref{tildeconv}) it will suffice to check that $X^N$ spends only an asymptotically negligible amount of time in $V(G^N)\backslash V(\tilde{G}^N)$. Since doing this requires only a simple adaptation of the proof of \cite[Lemma 3.8]{rwrrw}, we omit the details. To complete the proof, one then needs to replace $d_{{G}}$ by $d_{G^N}$, but this is straightforward since
\[N^{-1}\sup_{0\leq n\leq N}\left|d_{{G}}\left(0,S_n\right)-d_{G^N}\left(0,S_n\right)\right|\leq N^{-1}\left(T_1+T_{T_{N}^{-1}+1}-T_{T_{N}^{-1}}\right)\rightarrow0,\]
as $N\rightarrow\infty$.
\end{proof}

Although the previous lemma only contains a convergence statement for the random walks started from the particular vertex 0, there is no difficulty in extending this to the case when $X^N$ is started from a point $x_0^N\in V(G^N)$ such that $d_{G^N}(0,x_0^N)\sim \tau(d)^{-1}\delta(d)Nx_0$, and $B^{[0,1]}$ is started from $x_0\in[0,1]$.
Applying the local limit result of Proposition \ref{local} (to establish (\ref{tightcond}), we once again appeal to Lemma \ref{tightlem}), we are able deduce from this that Assumption \ref{assu1} holds for $\mathbf{P}$-a.e. realization of the original random walk.

{\lem For $\mathbf{P}$-a.e. realization of $S$, if $I\subset (0,\infty)$ is a compact interval, then
\[{\left(\left(V(G^N),\tau(d)\delta(d)^{-1}N^{-1}d_{G^N}\right), \pi^N, \left(q^N_{\kappa_2(d)^{-1}N^2t}(x,y)\right)_{x,y\in V(G^N),t\in I}\right)},\]
converges in $(\mathcal{M}_I,\Delta_I)$ to the triple consisting of: $[0,1]$ equipped with the Euclidean metric, Lebesgue measure on this set and the transition density of Brownian motion on $[0,1]$ reflected at the boundary.}
\bigskip

Since it is clear that \eqref{asum-a}-\eqref{erg} hold in this case, we can therefore apply Theorem \ref{main} to obtain the desired convergence of mixing times.

{\thm Fix $p\in[1,\infty]$. If $t_{\rm mix}^{p}(S_{[0,N]})$ is the  $L^p$-mixing time of the simple random walk on the range of $S$ up to time $N$, then $\mathbf{P}$-a.s.,
\[\kappa_2(d)N^{-2}t_{\rm mix}^{p}(S_{[0,N]}) {\rightarrow} {t_{\rm mix}^{p}([0,1])},\]
where $t_{\rm mix}^{p}([0,1])$ is the $L^p$-mixing time of the Brownian motion on $[0,1]$ reflected at the boundary.}

\section{Mixing time tail estimates}\label{tail section}

In this section, we give some sufficient conditions for deriving upper and lower estimates for mixing times of random walks on finite graphs, primarily using techniques adapted from \cite{NacPer08}. We will also discuss how to apply these general estimates to concrete random graphs (see Section \ref{exRG}). In order to crystallize the results and applications, most of the proofs shall be postponed to the appendix.

As will be illustrated by our examples, the results in this section are robust and convenient for obtaining mixing time tail estimates. Moreover, when the convergence of mixing times (as in Theorem \ref{main}) is available for a sequence of graphs, we highlight how, by first deriving estimates for the relevant continuous mixing time distribution (where similar techniques are sometimes applicable, see Remark \ref{upperrem}), it can be possible to deduce results regarding the asymptotic tail behavior of random graph mixing times that are difficult to obtain directly (see the proof of Proposition \ref{treecor} or Remark \ref{rem6-3ee}, for example).

\medskip

We start by fixing our notation. Let $G=(V(G), E(G))$ be a finite connected graph and $\mu^G$ be a weight function, as in the introduction. Suppose here that $d_G$ is the shortest path metric on the graph $G$, and denote, for a distinguished vertex $\rho\in V(G)$,
\[
     B(R)=\{y: d_G(\rho,y) <  R\}, ~
     V(R):=\sum_{x\in B(R)}\sum_{y:y\sim x} \mu^G_{xy} =\pi^G(B(R)) \mu(G),
     \quad R \in (0,\infty),
\]
where we write $x\sim y$ if $\mu_{xy}^G>0$ and set $\mu(G):=\sum_{x,y\in V(G)}\mu^G_{xy}$.
For the Markov chain $X^G$, let
\[
    \tau_R = \tau_{B(\rho,R)} =\min\{ n \ge 0: X^G_n \not\in B(R)\}.
\]
We define a quadratic form $\Ecal$ by
\[
      \Ecal(f,g)=\tfrac 12\sum_{\substack{x,y\in V(G) \\ x \sim y}}\mu_{xy}^G
      (f(x)-f(y))(g(x)-g(y)),
\]
and let $H^2=\{ f\in  \bR^{V(G)}: \Ecal(f,f)<\infty\}$. For disjoint subsets $A,B$ of $G$, the effective resistance between them is then given by:
\begin{equation}\label{eq:residef}
 \Reff(A,B)^{-1}=\inf\{\Ecal(f,f): f\in H^2, f|_A=1, f|_B=0\}.\end{equation}
If we further define $\Reff(x,y)=\Reff(\{x\},\{y\})$, and $\Reff(x,x)=0$, then one can check that $\Reff(\cdot,\cdot)$ is a metric on $V(G)$ (see
\cite[Section~2.3]{Kigami}). We will call this the resistance metric.
The resistance metric enjoys the following important (but easy to deduce) estimate,
\[|f(x)-f(y)|^2\le \Reff(x,y)\Ecal (f,f), \qquad\forall f\in L^2(G,\mu^G).\]
Moreover, it is easy to verify that if $c_1^{-1}:=\inf_{x,y\in G: x\sim y}\mu^G_{xy}$ $>0$, then
\begin{equation}\label{res-distre}
R_{\rm eff}(x,y)\le c_1d_G(x,y)\qquad \forall x,y\in G.
\end{equation}

Let $v, r: \{0,1,\cdots, {\rm diam}_{d_G}(G)+1\}\to [0,\infty)$ be strictly increasing
functions with $v(0)=r(0)=0$, $v(1)=r(1)=1$, which satisfy
\eq\label{volrescd}
C_1^{-1}\Big(\frac R{R'}\Big)^{d_1}\le \frac{v(R)}{v(R')}\le C_1
\Big(\frac R{R'}\Big)^{d_2},~~
C_2^{-1}\Big(\frac R{R'}\Big)^{\al_1}\le \frac{r(R)}{r(R')}\le C_2
\Big(\frac R{R'}\Big)^{\al_2}\en
for all $0<R'\le R\le {\rm diam}_{d_G}(G)+1$, where $C_1,C_2\ge 1$, $1\le d_1\le d_2$
and $0<\al_1\le \al_2\le 1$.
In what follows, $v(\cdot)$ will give the volume growth order and $r(\cdot)$ the resistance growth order.
For convenience, we extend them to functions on $[0,{\rm diam}_{d_G}(G)+1]$
by linear interpolation.
For the rest of the paper, $C_1,C_2, d_1,d_2$ and $\alpha_1,\alpha_2$ stand for the constants given in \eqref{volrescd}.

\subsection{General upper and lower bounds}
In this subsection, we give general upper and lower bounds for mixing times.
Note that, since $t^p_{\rm mix}(\rho)\le t^p_{\rm mix}(G)$ and $t^p_{\rm mix}(G)\le t^{p'}_{\rm mix}(G)$ for $p\le p'$,
it will be enough to estimate $t^\infty_{\rm mix}(G)$ for the upper bound\footnote{In fact, for the upper bound it is enough to estimate
$t^2_{\rm mix}(G)$. Indeed,
the Cauchy-Schwarz inequality
and \eqref{l2obs} implies the following known fact for mixing times of
symmetric Markov chains; $t^\infty_{\rm mix}(G)\le 2\,t^2_{\rm mix}(G; 1/2)$,
where $t^2_{\rm mix}(G; 1/2)$ is the $L^2$-mixing time of $G$ with $1/2$ instead of $1/4$ in the definition \eqref{mixdefdis}.}
and $t^1_{\rm mix}(\rho)$ for the lower bound.

\bigskip

\noindent
\underline{Upper bound}~
We first give an upper bound of the mixing times that is a reworking of \cite[Corollary 4.2]{NacPer08}, in our setting.
{\lem\label{upperCro} For any weighted graph $(G,\mu^G)$,
\[t^\infty_{\rm mix}(G)\leq 4{\rm diam}_{R}(G)\mu (G),\]
where ${\rm diam}_{R}(G)$ is the diameter of $G$ with respect to the resistance metric $\Reff$.}

\bigskip

\noindent
\underline{Lower bound}~
We next give the mixing time lower bound.
Let $\lam\ge 1, H_0,\cdots, H_3>0$, and let $C_3:=2^{-2/\al_1}C_2^{-1/\al_2}$ where $C_2$ is the constant in \eqref{volrescd}.
We give the following two conditions concerning the volume and resistance growth.
\eq\label{Cond1-1}
\Reff(\rho,y)\le \lam^{H_0}r(d_G(\rho,y)),~\forall y\in B(R),~~\mbox{ and }~~
V(R)\le \lam^{H_1}v(R),
\en
\eq\label{Cond1-2}
\Reff(\rho,B(R)^c)\ge \lam^{-H_2}r(R)~~\mbox{ and }~~
V(C_3\lam^{-(H_0+H_2)/\al_1}R)\ge \lam^{-H_3}v(C_3\lam^{-(H_0+H_2)/\al_1}R).
\en

{\propn \label{lowermn}
i) For $\lambda, R>1$, assume that $\mu(G)\ge 4V(R)$, and that
\eqref{Cond1-1}, \eqref{Cond1-2} hold for $R$, then
\eq
\label{mixlow-0}
t_{\rm mix}^1(G)>C_4\lam^{-H_2'-H_3}v(R) r(R),
\en
where $H_2'=H_2+(H_0+H_2)d_2/\al_1$. \\
ii) For $\lambda, R>1$, assume that $\mu(G)\ge 4V(R)$, and
\eqref{Cond1-1}, \eqref{Cond1-2} hold for $R$ and $\eps_0(\lam) R$, where
$\eps_0(\lam):=c_1\lam^{-(H_0+\sum_{i=0}^3H_i+H_2')/\al_1}$ for some $c_1>0$ small enough. Then
 \[t_{\rm mix}^1(\rho)>C_4\lam^{-H_2'-H_3}v(\eps_0(\lam) R) r(\eps_0(\lam) R).\]
}
{\rem \label{upperrem}{\rm
Essentially the same argument can be applied to deduce the corresponding mixing time upper and lower bounds in the continuous setting when we suppose that we have a process whose Dirichlet form is a resistance form.
(Remark \ref{upperrem2-0} contains details of the upper bound, and the details of the lower bound are omitted to avoid repetition.)

\subsection{Random graph case}

We now consider a probability space $(\Omega, \Fcal, \mathbf{P})$ carrying a
family of random weighted graphs $G^N(\omega)=(V(G^N(\omega)),E(G^N(\omega)), \mu^{N(\omega)}; \omega\in \Omega)$.
We assume
that, for each $N\in \mathbb{N}$ and $\omega \in \Omega$, $G^N(\omega)$ is a finite, connected graph containing a marked vertex $\rho^N$, and $\# V(G^N(\omega))\le M_N$ for some non-random constant $M_N<\infty$. (Here, for a set $A$, $\# A$ is the number of elements in $A$.)
Let $d_{G^N(\omega)}(\cdot,\cdot)$ be a graph distance, $B(R) := B_\omega(\rho^N,R)$, and $V(R) := V_\omega(\rho^N,R)$.
We write $X=(X_n , n \ge 0, P_\omega^x, x \in G^N(\omega))$ for the random walk on $G^N(\omega)$, and denote
by $p_n^\omega(x,y)$ its transition density with respect to $\pi^\omega$. Furthermore, we introduce a strictly increasing
function $h: \mathbb{N}\cup\{0\}\to [0,\infty)$ with $h(0)=0$, which will roughly describe the diameter of $G^N$
with respect to the graph distance. We then set $\gamma(\cdot)=v(h(\cdot))\cdot r(h(\cdot))$.
Finally, for $i=1,2$, we suppose $p_i: [1,\infty)\to [0,1]$ are functions such that $\lim_{\lam\to\infty}p_i(\lam)=0$. We then have the following.
(Note that $C_2, d_2$ in the statement are the constant in \eqref{volrescd}.)
{\propn\label{random-ul}
(1) Suppose that the following holds:
\eq\label{cond11-1}
{\mathbf P}({\rm diam}_R (G^N)\ge \lam r(h(N)))\le p_1(\lam),~~
{\mathbf P}(\mu^{N}(G^N) \ge \lam  v(h(N)))\le p_2(\lam),
\en
then
\[{\mathbf P}(t_{\rm mix}^\infty (G^N)\ge \lam \gamma(N))\le \inf_{\theta\in[0,1]} (p_1(\lam^{\theta}/8)+ p_2(\lam^{1-\theta})).\]
(2) Suppose there exist $c_1\le 1$ and $J\ge (1+H_1)/d_2$ such that the following holds:
\[{\mathbf P}(\mbox{\eqref{Cond1-1} $\wedge$ \eqref{Cond1-2} for } R=c_1\lam^{-J}h(N))\ge 1-p_1(\lam),~~
{\mathbf P}(\mu^N(G^N)< \lam^{-1}  v(h(N)))\le p_2(\lam),\]
then there exist $c_2,p_0>0$ such that
\[{\mathbf P}(t_{\rm mix}^1(G^N)\le c_2\lam^{-p_0}\gamma(N))\le 2p_1(\lam)+p_2(\lam/(4C_1c_1^{d_2})).\]
(3)  Suppose there exist $c_1\le 1$ and $J\ge (1+H_1)/d_2$ such that the following holds:
\begin{eqnarray}
\label{cond11-3-}
&&{\mathbf P}(\mbox{\eqref{Cond1-1} $\wedge$ \eqref{Cond1-2} for } R=c_1\lam^{-J}h(N) \mbox{ and for }
\eps_0(\lam)R)\ge 1-p_1(\lam),\nonumber\\
&&{\mathbf P}(\mu^N(G^N)< \lam^{-1}  v(h(N)))\le p_2(\lam),
\end{eqnarray}
where $\eps_0(\lam)$ is as in Proposition \ref{lowermn}\,ii),
then there exist $c_2,p_0>0$ such that
\[{\mathbf P}(t_{\rm mix}^1(\rho^N)\le c_2\lam^{-p_0}\gamma(N))\le 2p_1(\lam)+p_2(\lam/(4C_1c_1^{d_2})).\]
}

To illustrate this result, we consider the case when the random graphs $G^N(\omega)$ are obtained as components of percolation processes on finite graphs, thereby recovering \cite[Theorem 1.2(c)]{NacPer08}. (In \cite{NacPer08}, it was actually the lazy random walk was considered to avoid parity concerns, but the same techniques apply when we consider $q_m^G(\cdot,\cdot)$ as in \eqref{smooth} instead.)

{\propn \label{As-tau} Let $\hat G^N$
be a graph with $N$ vertices and with the maximum degree $d\in [3,N-1]$. Let ${\cal C}^N$ be the largest
component of the percolation subgraph of $\hat G^N$ for $0<p<1$. Let $p\leq \frac{1+\lambda n^{-1/3}}{d-1}$ for some
fixed $\lambda\in\mathbb{R}$, and assume that there exist $c_1,\theta_1\in(0,\infty)$ and $K_1\in \mathbb{N}$ such that
\begin{equation}\label{cdisnw}
{\mathbf P}(\#{\cal C}^N\leq A^{-1} N^{2/3}) \leq c_1A^{-\theta_1},\qquad\forall A, N\ge K_1,\end{equation}
then there exist $c_2, \theta_2\in(0,\infty)$ and $K_2\in\mathbb{N}$ such that, for all $p\in [1,\infty]$,
\begin{equation}\label{cdisnw-qq}
\mathbf{P}(A^{-1} N \leq t_{\rm mix}^{p}({\cal C}^N) \leq A N) \geq 1-c_2A^{-\theta_2},\qquad \forall  A, N\ge K_2.
\end{equation}}

\bigskip

Finally, below is a list of exponents for each example in Section \ref{examsec}.
\begin{center}
\begin{tabular}{|c|l|l|l|l|}
\hline
Section  & $v(R)$ & $r(R)$ & $h(N)$ & $\gamma (N)$\\
\hline
\ref{sssec} & $R^{\log K/\log L}$ & $R^{\log \lam/\log L}$ & $L^N$ &$(K\lam)^N$ \\
\hline
 \ref{treesec} with $a_N=N^{1/\alpha}, \alpha\in (1,2]$ & $R^{\alpha/(\alpha-1)}$ & $R$ & $N^{1-1/\alpha}$ & $N^{2-1/\alpha}$ \\
\hline
\ref{ersec} & $R^2$ & $R$ & $N^{1/3}$ & $N$ \\
\hline
\ref{rrwsec} & $R$  & $R$ & $N$ & $N^2$ \\
\hline
\end{tabular}
\end{center}
Here the Euclidean distance is used instead of the intrinsic shortest path metric for the examples in Section \ref{sssec}. Note that when $\alpha=2$ in Section \ref{treesec} (the finite variance case), the growth of $v(R)$ and $r(R)$ is of the same order as in Section \ref{ersec}. The difference of scaling exponents of mixing times (namely $\gamma (N)$) is due to the difference of scaling exponents for graph distances (namely $h(N)$). We also observe that the convergence to a stable law at \eqref{stable} forces the scaling constants to be of the form $a_N=N^{1/\alpha}L(N)$ for some slowly varying function $L$ (see \cite[Section 35]{GneKol}), and hence the above table captures all the most important first order behavior for the examples in Section \ref{treesec}.

\subsection{Examples}\label{exRG}

\noindent
\underline{Critical Galton-Watson trees of Section \ref{treesec}}~
By combining the results in this section with our mixing time convergence result, we can establish asymptotic bounds for the distributions of
mixing times of graphs in the sequence $(\mathcal{T}_N)_{N\geq 1}$ in the case when we have a finite variance offspring distribution.

{\propn\label{treecor} In the case when the offspring distribution has finite variance, there exist constants $c_1,c_2,c_3,c_4\in(0,\infty)$ such that
\begin{equation}\label{tupper}
\limsup_{N\rightarrow\infty} \mathbf{P}\left(N^{-3/2}t_{\rm mix}^{\infty}(\mathcal{T}_N)\geq\lambda\right)\leq c_1e^{-c_2\lambda^2}, \hspace{20pt}\forall\lambda\geq 0,
\end{equation}
and also
\begin{equation}\label{tlower}
\limsup_{N\rightarrow\infty}\mathbf{P}\left(N^{-3/2}t_{\rm mix}^{1}(\rho_N)\leq\lambda^{-1}\right)\leq c_3e^{-c_4\lambda^{1/25}}, \hspace{20pt}\forall\lambda\geq 0.
\end{equation}}
\begin{proof} To prove \eqref{tupper}, we apply the general mixing time upper bound of Lemma \ref{upperCro} to deduce that
\[\mathbf{P}\left(N^{-3/2}t_{\rm mix}^{\infty}(\mathcal{T}_N)\geq\lambda\right)\leq
\mathbf{P}\left(8N^{-1/2}{\rm diam}_{d_{\mathcal{T}_N}} (\mathcal{T}_N)\geq\lambda\right),\]
where ${\rm diam}_{d_{\mathcal{T}_N}} (\mathcal{T}_N)$ is the diameter of $\mathcal{T}_N$ with respect to $d_{\mathcal{T}_N}$, and we note that $\#E(\mathcal{T}_N)$ is equal to $2(N-1)$. By (\ref{treeconv}), the right-hand side here converges to $\mathbf{P}(8\,{\rm diam}_{d_{\mathcal{T}^{(2)}}} (\mathcal{T}^{(2)})\geq\lambda)$. By construction, the diameter of the continuum random tree $\mathcal{T}^{(2)}$ is bounded above by twice the supremum of the Brownian excursion of length 1. We can thus use the known distribution of the latter random variable (see \cite{Ken}, for example) to deduce the relevant bound.

For \eqref{tlower}, we first apply the convergence in distribution of Theorem \ref{treethm} to deduce that
\[\limsup_{N\rightarrow\infty}\mathbf{P}\left(N^{-3/2}t_{\rm mix}^{1}(\rho_N)\leq\lambda^{-1}\right)\leq \mathbf{P}\left(t_{\rm mix}^{1}(\rho)\leq\lambda^{-1}\right).\]
Now, for the continuum random tree, define
\[J(\lambda)=\{r>0 : \lambda^{-1}r^2\leq \pi^{(2)}(B_{\mathcal{T}^{(2)}}(\rho,r))\leq \lambda r^2,\,
R_\mathcal{T}^{(2)}(\rho,B_{\mathcal{T}^{(2)}}(\rho,r)^c)\geq \lambda^{-1}r\},\]
where $R_{\mathcal{T}^{(2)}}$ is the resistance on the continuum random tree (see \cite[(20)]{Croydoncrt}). Then
\[\mathbf{P}(r\in J(\lambda)) \geq 1-e^{-c\lambda},\hspace{20pt}\forall r\in(0,\tfrac{1}{2}],\lambda\geq 1,\]
(see \cite[Lemmas 4.1 and 7.1]{Croydoncrt}).
As a consequence of this, we can apply the continuous version of the mixing time lower bound discussed in Remark \ref{upperrem} (with $H_0=0$, $H_1=H_2=H_3=1$, $H_2'=3$, $\alpha_i=1$ and $d_i=2$) to deduce the desired result.
\end{proof}

{\rem {\rm The above proof already gives an estimate for the lower tail of $t_{\rm mix}^{1}(\rho)$. That the bound corresponding to (\ref{tupper}) holds for the limiting tree, i.e.
\[\mathbf{P}\left(t_{\rm mix}^{\infty}(\mathcal{T}^{(2)})\geq\lambda\right)\leq c_1e^{-c_2\lambda^2},\]
can be proved similarly to the discrete case (see Remark \ref{upperrem}).}}

\bigskip

\noindent
\underline{Critical Erd\H{o}s-R\'{e}nyi random graph of Section \ref{ersec}}~
Let $\mathcal{C}^N$ be the largest component of the Erd\H{o}s-R\'{e}nyi random graph in the critical window. Then the following holds.

{\propn\label{ERtaild}
There exist constants $c_1,c_2,c_3,N_0,\theta\in(0,\infty)$ such that
\begin{eqnarray}
\sup_{N\geq 1}\mathbf{P}\left(N^{-1}t_{\rm mix}^{\infty}(\mathcal{C}^N)\geq\lambda\right)&\leq & c_1e^{-c_2\lambda}, \hspace{20pt}\forall\lambda\geq 0,
\label{cupperfd}\\
\sup_{N\geq N_0}\mathbf{P}\left(N^{-1}t_{\rm mix}^{1}(\mathcal{C}^N)\leq\lambda^{-1}\right)&\leq & c_3\lambda^{-\theta}, \hspace{20pt}\forall\lambda\geq 0.
\label{clower}
\end{eqnarray}}
\begin{proof}
By \cite[Proposition 1.4]{NacPer08} and \cite[Theorem 1]{NacPer10}, \eqref{cupperfd} is an application of Proposition \ref{random-ul} with $p_1(\lambda)=c_4e^{-c_5\lambda^{3/2}}$ and $p_2(\lambda)=c_6e^{-c_7\lambda^{3}}$. \eqref{clower} is a consequence of Proposition \ref{As-tau}.
\end{proof}
{\rem\label{rem6-3ee} {\rm
(1) The tail estimates for $t_{\rm mix}^{1}(\mathcal{C}^N)$ are given in \cite[Theorem 1.1]{NacPer08} without quantitative bounds.
(In fact, reading the paper very carefully, it can be checked that the bounds similar to Proposition \ref{ERtaild} are available in the paper.) \\
(2) It does not seem possible to apply current estimates for the graphs $(\mathcal{C}^N)_{N\geq 1}$ and techniques for bounding mixing times to replace $t_{\rm mix}^{1}(\mathcal{C}^N)$ by $t_{\rm mix}^{1}(\rho^N)$
in the latter estimate (see Remark \ref{lremlwa}),
or even prove that the sequence $(N/t_{\rm mix}^{1}(\rho^N))_{N\geq 1}$ is tight, i.e.
\[\lim_{\lambda\rightarrow\infty}\limsup_{N\rightarrow\infty}\mathbf{P}\left(N^{-1}
t_{\rm mix}^{1}(\rho^N)\leq\lambda^{-1}\right)=0.\]
That this final statement is nonetheless true is a simple consequence of Theorem \ref{c1nmix}.
}}

\appendix
\section{Appendix: Proof of the statements in Section \ref{tail section}}
In this appendix, we prove various results given in Section \ref{tail section}.
We adopt the convention that if we cite elsewhere the constant $c_1$ of Proposition \ref{rw-keylemtau} (for example), we denote it as $c_{\ref{rw-keylemtau}.1}$.

\subsection{Proof of Lemma \ref{upperCro} and Proposition \ref{lowermn}}
~\indent {\it Proof of Lemma \ref{upperCro}.}~
First, note that by \cite[Proposition 3 in Chapter 2]{AF09}, we have that
\begin{equation}\label{ergo}\mathbf{E}^G_x\left(\sum_{m=0}^{\infty}\mathbf{1}_{\{X_m^G=x, m<S\}}\right)=\pi(x)\mathbf{E}^G_x(S),
\end{equation}
for any stopping time $S$ with $X^G_S=x$.
Taking $S$ to be the first hitting time of $x$ after time $2m-1$,
and writing $\Pi(x,2m)$ to represent the law of $X^G_{2m}$ when $X^G$ is started from $x$, we obtain that
\[\mathbf{E}^G_{\Pi(x,2m)}(\sigma_x)=\sum_{l=0}^{2m-1} \left(p^G_l(x,x)-1\right)
= 2\sum_{l=0}^{m-1}\left(q_{2l}^G(x,x)-1\right)\geq 2m\left( q^G_{2m}(x,x)-1\right),\]
where $\sigma_x$ is the first hitting time of $x$, and the inequality holds because $q^G_{2l}(x,x)$ is decreasing in $l$ (see the proof of \cite[Lemma 9]{CHLLT}, for example). Since by Cauchy-Schwarz, $|q^G_{2m}(x,y)-1|\leq (q^G_{2m}(x,x)-1)^{1/2}(q^G_{2m}(y,y)-1)^{1/2}$, it follows that
\[
\sup_{x\in V(G)}D_{\infty}^G(x,2m)=\sup_{x,y\in V(G)}\left|q_{2m}^G(x,y)-1\right|\leq  \sup_{x\in V(G)}(q^G_{2m}(x,x)-1)\leq \sup_{x,y\in V(G)}\frac{\mathbf{E}^G_{x}(\sigma_y)}{2m}.\]
By applying the commute time identity for random walks on graphs, $\mathbf{E}^G_{x}(\sigma_y)+\mathbf{E}^G_{y}(\sigma_x)=\Reff(x,y)\mu(G)$, this implies
$\sup_{x\in V(G)}D_{\infty}^G(x,2m)\leq {{\rm diam}_{R}(G) \mu(G)}/{2m}$, and the result follows.
\qed
{\rem \label{upperrem2-0}{\rm
As mentioned in Remark \ref{upperrem}, we can apply
essentially the same argument to deduce the corresponding mixing time upper bound in the continuous setting when we suppose that we have a process whose Dirichlet form is a resistance form. In particular, suppose that this is the case for $X^F$, as defined in the introduction. Let $S$ be the first hitting time of $x\in F$ after time $t$, then, for any $f\in L^1(F,\pi)$,
\begin{eqnarray*}
\mathbf{E}_{x}\left(\int_{0}^{S}f(X_s)ds\right)&=& \|f\|_{L^1(\pi)} \mathbf{E}_x(S),
\end{eqnarray*}
which can be obtained by applying an ergodicity argument similar to that used to prove (\ref{ergo}).  Writing $\Pi(x,t)$ to represent the law of $X^F_{t}$ when $X^F$ is started from $x$, the expectation on the right-hand side here satisfies
$\mathbf{E}_x(S)=t+\mathbf{E}_{\Pi(x,t)}(\tau_x)\leq t+\sup_{y\in F}\Reff(x,y)$,
where to deduce the upper bound, we have applied that the commute time identity $\mathbf{E}_{x}(\tau_y)+\mathbf{E}_{y}(\tau_x)=\Reff(x,y)$ also holds for resistance forms (since we are assuming $\pi$ to be a probability measure, it does not appear explicitly in this version of the identity). Moreover, if $f$ is positive, the left-hand side is bounded below as follows: $\mathbf{E}_{x}(\int_{0}^{S}f(X_s)ds)\geq \int_{0}^t\int_F q_s(x,y)f(y)\pi(dy)ds$. Combining these bounds, we have proved that, for positive $f\in L^1(F,\pi)$ such that $ \|f\|_{L^1(\pi)} \neq 0$,
\[\frac{\int_{0}^t\int_F q_s(x,y)f(y)\pi(dy)ds}{\|f\|_{L^1(\pi)}}\leq t+{\rm diam}_{R}(F).\]
By choosing a sequence of suitable functions whose support converges to $\{x\}$, the joint continuity of $(q_t(x,y))_{x,y\in F,t>0}$ allows us to deduce from this that
\[tq_t(x,x)\leq\int_0^t q_s(x,x)ds\leq t+{\rm diam}_{R}(F),\]
where the first inequality holds because $q_t(x,x)$ is decreasing in $t$. The remainder of the proof
is identical to the graph case.
}
}

\bigskip

The proof of Proposition \ref{lowermn} requires some preparations.
Our argument depends on some estimates for hitting times that are modifications of results in \cite{bjks,KumMi}.

To begin with, let $B=B(R)$ and define
\[
  g_B(x,y) = \mu_y^{-1} \sum_{k=0}^\infty \mathbf{P}^G_x( X_k=y, k<\tau_B).
\]
Then, it is easy to show that
\[\mathbf{E}^G_z\tau_B = \sum_{y\in B} g_B(z,y)\mu_y,\qquad
\Reff(x,B^c)=g_B(x,x)\]
(see, for example \cite[(2.19),(2.20)]{bjks}).
Also, if $A$ and $B$ are disjoint subsets of $G$ and $x \notin A \cup B$, then
(see \cite[(2.14)]{bjks})
\eq \label{bgp-lem}
\mathbf{P}^G_x(T_A < T_B) \le \frac{\Reff(x,B)}{\Reff(x,A)},
\en
where $T_A$ is the hitting time of $A \subset G$. If $C_4:=8^{-1}C_1^{-1}C_3^{d_2}$, we can then prove the following.
(Here, recall that  $C_3=2^{-2/\al_1}C_2^{-1/\al_2}$ and $C_1, C_2$ are the constants in \eqref{volrescd}.)

{\lem
Let $\lam\ge 1$ and $H_0,\cdots, H_3>0$.\\
(a) Suppose \eqref{Cond1-1} holds. Then
\eq
\label{Etaub}
    \mathbf{E}^G_x \tau_R \le 2 \lam^{H_0+H_1} v(R)r(R) \quad\text{ for } x\in B(R).\\
\en
(b) Suppose \eqref{Cond1-1} and \eqref{Cond1-2} hold.
Then
\begin{align}
\label{e:tmlb}
   \mathbf{E}^G_x  \tau_R &\ge 2C_4\lam^{-H_2'-H_3}v(R) r(R)~~ \text { for } x \in  B(C_3\lam^{-(H_0+H_2)/\al_1}R),
\end{align}
where we recall $H_2'=H_2+(H_0+H_2)d_2/\al_1$. \\
(c) Suppose \eqref{Cond1-1} and \eqref{Cond1-2}, and let $x \in B(C_3\lam^{-(H_0+H_2)/\al_1}R)$, then
\eq \label{e:tplb}
    \mathbf{P}^G_x ( \tau_R > n) \ge \frac{2C_4\lam^{-H_2'-H_3}v(R) r(R) -n}{2\lam^{H_0+H_1}v(R) r(R)}  \quad
\text{ for } n \ge 0. \\
\en}
\begin{proof}
Using \eqref{Cond1-1}, we have
$\Reff(z,B^c)\le \Reff(0,z)+\Reff(0,B^c)\le 2\lam^{H_0}r(R)$ for any
$z\in B$. So,
\[ \mathbf{E}^G_z\tau_B= \sum_{y\in B} g_B(z,y)\mu_y \le \sum_{y\in B} g_B(z,z) \mu_y
  = \Reff(z,B^c) V(R) \le 2\lam^{H_0+H_1} v(R)r(R),\]
{which gives (\ref{Etaub})}. In order to prove \eqref{e:tmlb}, we first establish the
following: for $0< \eps \le 1/(2C_2\lam^{H_0+H_2})^{1/\al_1}$
$=2^{1/\al_1}C_3\lam^{-(H_0+H_2)/\al_1}$ and $y \in B(\eps R)$,
we have
\begin{align}\label{e:dowq}
 \mathbf{E}^G_y ( T_\rho < \tau_R ) &\ge 1-  \frac{C_2\eps^{\al_1}\lam^{H_0+H_2}}{1-C_2\eps^{\al_1}\lam^{H_0+H_2}}
\ge 1- 2 C_2\eps^{\al_1}\lam^{H_0+H_2}.
\end{align}
Indeed, by the first inequalities of \eqref{Cond1-1} and \eqref{Cond1-2},
we have
\[ \Reff(y,B(R)^c) \ge \Reff(\rho,B(R)^c) - \Reff(\rho,y)
 \ge \lam^{-H_2}r(R) - \lam^{H_0}r(\eps R) \ge \frac{r(\eps R)}{C_2\eps^{\al_1}\lam^{H_2}} - \lam^{H_0}r(\eps R). \]
So, by \eqref{bgp-lem},
\[ \mathbf{P}^G_y(\tau_R < T_\rho)
\le \frac{\Reff(y,\rho)}{\Reff(y, B(R)^c)}
\le  \frac { \lam^{H_0}r(\eps R)}{\frac{r(\eps R)}{C_2\eps^{\al_1}\lam^{H_2}} - \lam^{H_0}r(\eps R)}
\le \frac{C_2\eps^{\al_1}\lam^{H_0+H_2}}{1-C_2\eps^{\al_1}\lam^{H_0+H_2}},\]
and \eqref{e:dowq} is obtained.
Now, if $y \in B'=B(C_3\lam^{-(H_0+H_2)/\al_1}R)$, then the bound at \eqref{e:dowq} gives that
$ \mathbf{P}^G_y(T_\rho < \tau_B) \ge {\tfrac12}$, so
\[ g_B(\rho,y) = g_B(\rho,\rho) \mathbf{P}^G_y (T_\rho < \tau_B) \ge {\tfrac12}  g_B(\rho,\rho)
 = \tfrac12 \Reff(\rho,B^c) \ge \tfrac 12  \lam^{-H_2}r(R). \]
By the second inequality of \eqref{Cond1-2}, we have
\[\mu(B') \ge \lam^{-H_3}v(C_3\lam^{-(H_0+H_2)/\al_1}R)\ge
C_1^{-1}C_3^{d_2}\lam^{-2(H_0+H_2)d_2/\al_1-H_3}v(R),\]
and therefore we obtain,
\begin{eqnarray}
  \mathbf{E}^G_\rho\tau_B &\ge& \sum_{y\in B'} g_B(\rho,y)\mu_y\nonumber\\
 &\ge& {\tfrac12} g_B(\rho,\rho) \mu(B')\nonumber\\
 &\ge &\tfrac 12
C_1^{-1}C_3^{d_2}\lam^{-H_2-(H_0+H_2)d_2/\al_1-H_3}v(R) r(R)\nonumber\\
&=&4C_4\lam^{-H_2'-H_3}v(R) r(R).\nonumber
\end{eqnarray}
Moreover, for $x\in B'$ we have that $ \mathbf{E}^G_x  \tau_B \ge  \mathbf{P}^G_x (T_\rho < \tau_B)  \mathbf{E}^G_\rho \tau_B$, which gives
\eqref{e:tmlb}.

Finally, by the Markov property,
\eqref{Etaub} and \eqref{e:tmlb},
\begin{eqnarray*}
 2C_4\lam^{-H_2'-H_3}v(R) r(R)\le
 \mathbf{E}^G_x \tau_R &\le& n+  \mathbf{E}^G_x [\mathbf{1}_{\{\tau_R> n\} }  \mathbf{E}^G_{X_n}(\tau_R)]\\
&\le& n +
 2 \lam^{H_0+H_1} v(R)r(R)
  \mathbf{P}^G_x  ( \tau_R>n).
 \end{eqnarray*}
Rearranging this gives
\eqref{e:tplb}. \end{proof}

\smallskip

The following estimate is a modification of \cite[Proposition 3.5\,(a)]{KumMi} (see \cite[(2.4)]{bjks} for the important
special case $v(R)=R^2$, $r(R)=R$). Note that for $R>{\rm diam}_{d_G} (G)$, it is the case that $\tau_R=\infty$, and so \eqref{trdlb} trivially holds.

{\propn \label{rw-keylemtau}
Let $0< \eps \le C_3\lam^{-(H_0+H_2)/\al_1}$, and suppose \eqref{Cond1-1} and \eqref{Cond1-2} for $R$ and $\eps R$, then
\eq
 \label{trdlb}
  \mathbf{P}^G_y\big( \tau_R \le C_4\lam^{-H_2'-H_3}v(\eps R) r(\eps R) \big) \le c_1\lam^{H_0+\sum_{i=0}^3H_i+H_2'} \eps^{\al_1},
 \quad \text{ for } y \in B(\eps R).
\en
}

\begin{proof}
We take a kind of bootstrap from \eqref{e:tplb}
and \eqref{e:dowq}. Let $t_0>0$, and set
$$ q(y)= \mathbf{P}^G_y (\tau_R {\leq} T_\rho), \qquad a(y)=  \mathbf{P}^G_y (\tau_R \le t_0). $$
Then
\begin{align}
 a(y) =  \mathbf{P}^G_y(\tau_R \le t_0)
 &=  \mathbf{P}^G_y( \tau_R \le t_0, \tau_R {\leq} T_\rho)
 +  \mathbf{P}^G_y( \tau_R \le t_0,  \tau_R > T_\rho)\nonumber\\
 &\le   \mathbf{P}^G_y( \tau_R \le T_\rho) +   \mathbf{P}^G_y(T_\rho< \tau_R,  \tau_R -T_\rho \le t_0)\nonumber\\
\label{pytub}
 &\le q(y) + (1-q(y)) a(\rho)  \le q(y) + a(\rho),
\end{align}
using the strong Markov property for the second inequality.
Starting the Markov chain $X$ at $\rho$, we have
\begin{align}\label{pytub2}
a(\rho) =  \mathbf{P}^G_\rho(\tau_R\le t_0) \le
 \mathbf{E}^G_\rho [1_{\{\tau_{\eps R}\le t_0 \}}  \mathbf{P}^G_{X_{\tau_{\eps R}}}( \tau_R \le t_0)]
  \le   \mathbf{P}^G_\rho( \tau_{\eps R} \le t_0)
\max_{y \in  \partial B(\eps R )} a(y).
\end{align}
Combining  \eqref{pytub} and \eqref{pytub2} gives
\eq \label{a0ub}
  a(\rho) \le \frac{  \max_{y   \in \partial B(\eps R )} q(y)}
 {   \mathbf{P}^G_\rho( \tau_{\eps R} > t_0)}.
\en
Further, using \eqref{e:dowq} with $2\eps$, we have
\eq
\label{qy}
 q(y) \le \frac{C_2(2\eps)^{\al_1}\lam^{H_0+H_2}}{1-C_2(2\eps)^{\al_1}\lam^{H_0+H_2}}
 \le 2C_2(2\eps)^{\al_1}\lam^{H_0+H_2}.
\en
Let
$t_0=  C_4\lam^{-H_2'-H_3}v(\eps R) r(\eps R)$; then using \eqref{e:tplb} for the
ball $B(\eps R)$ (note that \eqref{Cond1-1} and \eqref{Cond1-2} for $\eps R$ are assumed here),
we obtain
$$  \mathbf{P}^G_\rho( \tau_{\eps R} > t_0) \ge c_0\lam^{-H_0-H_1-H_2'-H_3}. $$
combining this with \eqref{qy}, (\ref{a0ub}) and (\ref{pytub}) completes
the proof of \eqref{trdlb}. \end{proof}

Note that, we may and will take $c_{\ref{rw-keylemtau}.1}>1/(2C_3^{\al_1})$.
Now we are ready to prove Proposition \ref{lowermn}.

{\it Proof of Proposition \ref{lowermn}.}~
i) We follow the argument in \cite[Lemma 5.4]{NacPer08}.
Let $t\in \mathbb{N}$. If $ \mathbf{P}^G_x(\tau_B\le t)\ge 1/2$ for all $x\in B(R-1)$, then $\tau_R/t$ is stochastically
dominated by a geometric random variable with parameter $1/2$, so that $ \mathbf{E}^G_\rho[\tau_R]\le 2t$. By this and \eqref{e:tmlb},
we see that for $t=C_4\lam^{-H_2'-H_3}v(R) r(R)$, there exists some $x\in B(R-1)$ such that $ \mathbf{P}^G_x(\tau_B\le t)\le 1/2$.
Further, since  $\mu(G)\ge 4V(R)$,
$\pi(B(R))=V(B(R))/\mu(G)\le 1/4$.
Combining these observations, we obtain
\eq\label{eq:henql}
D_1(x,t)\geq2\mathbf{P}^G_x\left(\tau_{R}\geq t\right)-2\pi(B(R))\geq 1-\frac 12> \frac 14,\en
so that \eqref{mixlow-0} follows.

ii) Take $\eps=\eps_0(\lam)$ in Proposition \ref{rw-keylemtau} and let
 $t=C_4\lam^{-H_2'-H_3}v(\eps R) r(\eps R)$. Then, since
$0< \eps \le C_3\lam^{-(H_0+H_2)/\al_1}$ (this is because we take $c_{\ref{rw-keylemtau}.1}>1/(2C_3^{\al_1})$), by
\eqref{trdlb} we have ${\mathbf P}^G_\rho \big( \tau_R \le t \big) \le c_{\ref{rw-keylemtau}.1}\lam^{H_0+\sum_{i=0}^3H_i+H_2'} \eps^{\al_1}=1/2$.
The rest is the same as the proof of i) except that we take $x=\rho$ in \eqref{eq:henql} and take
$c_{\ref{lowermn}.1}=(2c_{\ref{rw-keylemtau}.1})^{-1/\al_1}$. \qed

\subsection{ Proof of Proposition \ref{random-ul} and Proposition \ref{As-tau}}
~\indent {\it Proof of Proposition \ref{random-ul}.}~
By Lemma \ref{upperCro}, we have for any $\theta\in[0,1]$ that
\begin{eqnarray*}
\mathbf{P}\left(t_{\rm mix}^{\infty}(G^N)\geq \lambda \gamma(N)\right)&\leq&
\mathbf{P}\left(8{\rm diam}_{R}(G^N)\mu^N(G^N)\geq \lambda \gamma(N)\right)\\
&\leq & \mathbf{P}\left(8{\rm diam}_{R}(G^N)\geq \lambda^{\theta} r(h(N))\right)+\mathbf{P}\left(\mu^N(G^N)\geq \lambda^{1-\theta}
v(h(N))\right)\\
&\leq & p_1(\lam^{\theta}/8)+ p_2(\lam^{1-\theta}),
\end{eqnarray*}
which implies the conclusion of (1).

For (2), let $R=c_1\lam^{-J}h(N)$ and define
\begin{eqnarray*}
t&:=&C_4\lam^{-H_2'-H_3}v(R) r(R)=C_4\lam^{-H_2'-H_3}v(c_1\lam^{-J}h(N)) r(c_1\lam^{-J}h(N))\\
&\ge & C_4\lam^{-H_2'-H_3}C_1^{-1}C_2^{-2}(c_1\lam^{-J})^{d_2+\al_2}v(h(N)) r(h(N))=:c_2\lam^{-p_0}\gamma(N).
\end{eqnarray*}
Then by Proposition \ref{lowermn}\,i),
\begin{eqnarray*}
\lefteqn{{\mathbf P}(t_{\rm mix}^1(G^N)\le c_1\lam^{-p_0}\gamma(N))\le {\mathbf P}(t_{\rm mix}^1(G^N)\le t)}\\
&\le &
{\mathbf P}(\mbox{either \eqref{Cond1-1} or \eqref{Cond1-2} do not hold for } R=c_1\lam^{-J}h(N))+{\mathbf P}(\mu^N(G^N)<4V(R))\\
&\le &p_1(\lam)+{\mathbf P}(\mu^N(G^N)<4V(R)).\end{eqnarray*}
Note that
\[4\lam^{H_1}v(R)=4\lam^{H_1}v(c_1\lam^{-J}h(N))\le 4\lam^{H_1}C_1(c_1\lam^{-J})^{d_2}v(h(N))\le
4C_1c_1^{d_2}\lam^{-1}v(h(N)),\]
where we used $J\ge (1+H_1)/d_2$ in the last inequality. Using this, we have
\begin{eqnarray*}
\lefteqn{{\mathbf P}(\mu^N(G^N)<4V(R))}\\&\le&
{\mathbf P}(\mu^N(G^N)<4\lam^{H_1}v(R))+{\mathbf P}(\lam^{H_1}v(R)\le V(R))\\
&\le & {\mathbf P}(\mu^N(G^N)<4C_1c_1^{d_2}\lam^{-1}v(h(N)))+p_1(\lam)\\
&\le& p_2(\lam/(4C_1c_1^{d_2}))+p_1(\lam), ~~\mbox{  }~ \end{eqnarray*}
which implies the conclusion of (2).
The proof of (3) is almost the same, so we omit it.\qed

{\it Proof of Proposition \ref{As-tau}.}~
We only indicate how to apply previous propositions.
First, the upper bound of $t_{\rm mix}^{p}({\cal C}^N)$ can be obtained by Proposition \ref{random-ul}\,(1)
with $v(R)=R^2, r(R)=R, h(N)=N^{1/3}$ and $p_1(A)=c_0A^{-q_0}, p_2(A)=c_0'A^{-q_0'}$
for some $c_0,c_0',q_0,q_0'>0$.
Indeed, \eqref{cond11-1} holds because of \cite[Theorem 2.1 (a),(b), Theorem 6.1]{NacPer08} and the fact
${\rm diam}\,({\cal C}^N)\ge {\rm diam}_R\,({\cal C}^N)$, which is due to \eqref{res-distre}.

The lower bound is more complicated.
Using Proposition 5.5--5.7 and (5.1) in \cite{NacPer08} with
\[\beta=\lam^{-1/4}, L=\lam^{H_2}, \al=\lam^{H_1}, r=R, h=C_3\lam^{-H_2}R, m=\lam^{-H_3}(C_3\lam^{-H_2}R)^2,\]
and then taking $R=c_1\lam^{-J}N^{1/3}$, $H_0=0$ (due to \eqref{res-distre}),
$H_1=H_2=2, H_3=4, J=(1+H_1)/2=3/2$, we see that for each $v\in \hat G^N$,
\[{\mathbf P}(\#{\cal C}(v)>\lam^{-1/4} N^{2/3} \mbox{ and } {\cal A})\le c_4\lam^{-1/2}N^{-1/3},\]
where
\[{\cal A}=
\{V(v,C_3\lam^{-2}R)\le \lam^{-5}(C_3\lam^{-2}R)^2,\,\Reff(v,B(v,R)^c)\le \frac R{8\lam^2},\, \#E(B(v,R))\ge \lam^2 R^2\}.\]
This corresponds to \cite[(5.3)]{NacPer08}. Now using Proposition \ref{lowermn}\,i) and arguing similarly to the
proof of \cite[Theorem 2.1 (c.2)]{NacPer08}, we have
\[{\mathbf P}(\exists v\in \hat G^N\mbox{ with }\#{\cal C}(v)>\lam^{-1/4} N^{2/3} \mbox{ and }
t_{\rm mix}^1({\cal C}(v))\le C_4\lam^{-29/2}N)\le  c_4\lam^{-1/4}.\]
This together with \eqref{cdisnw} implies the desired lower bound of $t_{\rm mix}^{p}({\cal C}^N)$.
\qed

The proofs of this proposition and Proposition \ref{treecor} highlight why it is useful to have a general theory where the exponents $H_0, \cdots, H_3$ can vary.

{\rem\label{lremlwa}{\rm
As mentioned in Remark \ref{rem6-3ee}\,(2),
it does not seem possible to apply current estimates for the graphs $(\mathcal{C}^N)_{N\geq 1}$ and techniques for bounding mixing times to replace $A^{-1} N \leq t_{\rm mix}^{p}(\mathcal{C}^N)$ by $A^{-1} N \leq t_{\rm mix}^{p}(\rho^N)$ in \eqref{cdisnw-qq}. The major difficulty
is to verify the first inequality of \eqref{cond11-3-} for $\eps_0(\lam)R$. Indeed, even if we choose $H_0,\cdots, H_3$ large
(which increases the chance that \eqref{Cond1-1} and \eqref{Cond1-2} hold for $R$), $\eps_0(\lam)$ gets small accordingly, so that
the probability ${\mathbf P}(\mbox{\eqref{Cond1-1} $\wedge$ \eqref{Cond1-2} for }\eps_0(\lam)R)$ does not increase.}}

\providecommand{\bysame}{\leavevmode\hbox to3em{\hrulefill}\thinspace}
\providecommand{\MR}{\relax\ifhmode\unskip\space\fi MR }
\providecommand{\MRhref}[2]{
  \href{http://www.ams.org/mathscinet-getitem?mr=#1}{#2}
}
\providecommand{\href}[2]{#2}


\begin{thebibliography}{10}

\bibitem{ABG}
L.~Addario-Berry, N.~Broutin, and C.~Goldschmidt, \emph{The continuum limit of
  critical random graphs}, Probab. Theory Related Fields, to appear.

\bibitem{AF09}
D. Aldous and J. Fill,
\emph{Reversible Markov chains and random walks on graphs},
\newblock Preprint
{\tt http://www.stat.berkeley.edu/$\sim$aldous/RWG/book.html}

\bibitem{bjks}
M.T. Barlow, A.A. J\'{a}rai, T. Kumagai and G. Slade,
\emph{Random walk on the incipient infinite cluster for oriented percolation
in high dimensions},
Comm. Math. Phys. \textbf{278} (2008), 385--431.

\bibitem{BKW}
I. Benjamini, G. Kozma and N. Wormald, \emph{The mixing time of the
giant component of a random graph}, preprint.

\bibitem{BBG}
P. B\'erard, G. Besson and S. Gallot,
 \emph{Embedding Riemannian manifolds by their heat kernel},
Geom. Funct. Anal. \textbf{4} (1994), 373--398.

\bibitem{BCHSS}
C. Borgs, J.T. Chayes, R. van der Hofstad, G. Slade and J. Spencer, Random subgraphs of
finite graphs: I. The scaling window under the triangle condition, \emph{Random Structures Algorithms,}
{\bf 27} 137--184, 2005.

\bibitem{BBI}
D.~Burago, Y.~Burago, and S.~Ivanov, \emph{A course in metric geometry},
  Graduate Studies in Mathematics, vol.~33, American Mathematical Society,
  Providence, RI, 2001.

\bibitem{S-C}
G.-Y. Chen and L.~Saloff-Coste, \emph{The cutoff phenomenon for ergodic
  {M}arkov processes}, Electron. J. Probab. \textbf{13} (2008), 26--78.

\bibitem{Croydoncrg}
D.~A. Croydon, \emph{Scaling limit for the random walk on the largest connected
  component of the critical random graph},
  Publ. RIMS. Kyoto Univ., to appear.

\bibitem{Croydoncbp}
\bysame, \emph{Convergence of simple random walks on random discrete trees to
  {B}rownian motion on the continuum random tree}, Ann. Inst. Henri Poincar\'e
  Probab. Stat. \textbf{44} (2008), 987--1019.

\bibitem{Croydoncrt}
\bysame, \emph{Volume growth and heat kernel estimates for the continuum random tree}, Probab. Theory Related Fields \textbf{140} (2008), 207--238.

\bibitem{rwrrw}
\bysame, \emph{Random walk on the range of random walk}, J. Stat. Phys.
  \textbf{136} (2009), 349--372.

\bibitem{Croydoninf}
\bysame, \emph{Scaling limits for simple random walks on random ordered graph
  trees}, Adv. in Appl. Probab. \textbf{42} (2010), 528--558.

\bibitem{CHLLT}
D.~A. Croydon and B.~M. Hambly, \emph{Local limit theorems for sequences of
  simple random walks on graphs}, Potential Anal. \textbf{29} (2008),
  351--389.

\bibitem{Duqap}
T.~Duquesne, \emph{A limit theorem for the contour process of conditioned
  {G}alton-{W}atson trees}, Ann. Probab. \textbf{31} (2003), 996--1027.

\bibitem{LegallDuquesne}
T.~Duquesne and J.-F. Le~Gall, \emph{Probabilistic and fractal aspects of
  {L}\'evy trees}, Probab. Theory Related Fields \textbf{131} (2005),
  553--603.

\bibitem{ET}
P.~Erd{\H{o}}s and S.~J. Taylor, \emph{Some intersection properties of random
  walk paths}, Acta Math. Acad. Sci. Hungar. \textbf{11} (1960), 231--248.

\bibitem{FR}
N. Fountoulakis and B.A. Reed,  \emph{The evolution of the mixing rate of a simple
random walk on the giant component of a random graph}, Random Structures Algorithms,
\textbf{33} (2008), 68--86.

\bibitem{FOT}
M. Fukushima, Y. Oshima and M. Takeda, \emph{Dirichlet forms and symmetric Markov processes},
de Gruyter Studies in Mathematics, 19. Walter de Gruyter \& Co., Berlin, 2011.

\bibitem{GneKol}
B.~V. Gnedenko and A.~N. Kolmogorov, \emph{Limit distributions for sums of
  independent random variables}, Addison-Wesley Publishing Company, Inc.,
  Cambridge, Mass., 1954, Translated and annotated by K. L. Chung. With an
  Appendix by J. L. Doob.

\bibitem{GMT}
S.~Goel, R.~Montenegro, and P.~Tetali, \emph{Mixing time bounds via the
  spectral profile}, Electron. J. Probab. \textbf{11} (2006), 1--26.

\bibitem{GPW}
A.~Greven, P.~Pfaffelhuber, and A.~Winter, \emph{Convergence in distribution of
  random metric measure spaces ({$\Lambda$}-coalescent measure trees)}, Probab.
  Theory Related Fields \textbf{145} (2009), 285--322.


\bibitem{HeyvdH09}
M. Heydenreich and R. van der Hofstad,  \emph{Random graph asymptotics on high-dimensional tori II. Volume, diameter and mixing time},
Probab. Theory Related Fields, \textbf{149} (2011), 397--415.

\bibitem{KasK}
A. Kasue and H. Kumura, \emph{Spectral convergence of Riemannian manifolds}, T\^ohoku Math. J.
\textbf{46} (1994), 147--179.

\bibitem{Ken}
D. P. Kennedy, \emph{The distribution of the maximum Brownian excursion}, J. Appl. Probab.
\textbf{13} (1976), 371--376.

\bibitem{KigRes}
J.~Kigami, \emph{Resistance forms, quasisymmetric maps and heat kernel
  estimates}, Memoirs AMS, to appear.

\bibitem{Kigamimetric}
J.~Kigami, \emph{Hausdorff dimensions of self-similar sets and shortest path metrics}, J. Math. Soc. Japan \textbf{47} (1995), 381--404.

\bibitem{Kigamidendrite}
J.~Kigami, \emph{Harmonic calculus on limits of networks and its application to
  dendrites}, J. Funct. Anal. \textbf{128} (1995), 48--86.

\bibitem{Kigami}
J.~Kigami, \emph{Analysis on fractals}, Cambridge Tracts in Mathematics, vol.
  143, Cambridge University Press, Cambridge, 2001.

\bibitem{homog}
T.~Kumagai, \emph{Homogenization on finitely ramified fractals}, Stochastic
  analysis and related topics in Kyoto, Adv. Stud. Pure Math., vol.~41, Math.
  Soc. Japan, Tokyo, 2004, pp.~189--207.

\bibitem{KumKus}
T. Kumagai and S. Kusuoka, \emph{Homogenization on nested fractals},
Probab. Theory Related Fields \textbf{104} (1996), 375-398.

\bibitem{KumMi}
T. Kumagai and J. Misumi,
\emph{Heat kernel estimates for strongly recurrent random walk on random media},
J. Theoret. Probab. {\bf 21} (2008), 910--935.

\bibitem{rrt}
J.-F. Le~Gall, \emph{Random real trees}, Ann. Fac. Sci. Toulouse Math. (6)
\textbf{15} (2006), 35--62.

\bibitem{LPW}
D. Levin, Y. Peres and E. Wilmer,
\emph{Markov chains and mixing times},
Amer. Math. Soc., Providence, RI, 2009.

\bibitem{NacPer08}
A. Nachmias and Y. Peres,
\emph{Critical random graphs: diameter and mixing time},
Ann. Probab. {\bf 36} (2008), 1267--1286.

\bibitem{NacPer10}
A. Nachmias and Y. Peres,
\emph{The critical random graph, with martingales},
Israel J. Math. 176 (2010), 29--41.

\end{thebibliography}
\end{document}